\documentclass[11pt]{amsart}

\usepackage{amsmath}
\usepackage{amsthm}
\usepackage{amsfonts}
\usepackage{amssymb}
\usepackage{verbatim}
\usepackage{graphicx}
\usepackage[margin=1.0in]{geometry}

\newcommand\R{{\mathbb{R}}}
\newcommand\C{{\mathbb{C}}}

\newcommand\I{{\mathbf{I}}}
\renewcommand\P{{\mathbf{P}}}
\newcommand\E{{\mathbf{E}}}

\renewcommand\Im{{\operatorname{Im}}}
\renewcommand\Re{{\operatorname{Re}}}
\newcommand\eps{{\varepsilon}}
\newcommand\trace{\operatorname{trace}}
\newcommand\tr{\operatorname{trace}}

\renewcommand\Pr{{\mathbf P }}

\newcommand\CP{{\mathcal P}}

\newcommand\BBC{{\mathbb C}}

\newcommand\ep{{\epsilon}}

\theoremstyle{plain}
  \newtheorem{theorem}{Theorem}[section]
  
  \newtheorem{problem}[theorem]{Problem}

  \newtheorem{lemma}[theorem]{Lemma}
  \newtheorem{corollary}[theorem]{Corollary}

\theoremstyle{definition}

  \newtheorem{remark}[theorem]{Remark}

\include{psfig}

\begin{document}

\title{Random weighted  projections, Random quadratic forms and Random  Eigenvectors}

\author[V. Vu and K. Wang]{Van Vu\\Department of Mathematics, Yale University\\ \\Ke Wang\\Institute for Mathematics and its Applications, University of Minnesota}
\thanks{V. Vu is supported by research grants DMS-0901216, DMS-1307797 and AFOSAR-FA-9550-09-1-0167. }
\address{Van Vu\\Department of Mathematics, Yale University, New Haven, CT 06520, USA}
\email{van.vu@yale.edu}

\address{Ke Wang\\Institute for Mathematics and its Applications, University of Minnesota, Minneapolis, MN  55455, USA }
\email{wangk@umn.edu}


\begin{abstract}
We present a concentration result concerning random weighted projections in high dimensional spaces. As applications, we prove

\begin{itemize} 

\item New  concentration inequalities for random quadratic forms.

\item  The infinity norm of most unit eigenvectors of  a random $\pm 1$ matrix is of  order $O( \sqrt { \log n/n})$.

\item An  estimate on the threshold for the local semi-circle law which is tight up to a $\sqrt {\log n}$ factor. \end{itemize} 
\end{abstract}


\maketitle

\maketitle

\section{Introduction}

\subsection{Projection of a random vector} 

Consider $\C^n$  with a subspace $H$ of dimension $d$. Let $X =(\xi_1, \dots, \xi_n)$ be a random vector. In all considerations in this paper, we assume that $X$ is in isotropic position, namely $\E X \otimes X= Id$. 
The length of the  orthogonal projection of $X$ onto 
$H$ is an important parameter which plays an essential role in the studies of random matrices and related areas.

In \cite{TVdet}, Tao and the first author showed that (under certain conditions) this length is strongly concentrated. In other words,   the projection 
of $X$ onto $H$ lies essentially on a circle centered at the origin. 
 This fact played a crucial role in the studies  of the determinant of a random matrix with independent entries (see \cite{TVdet, NVdet}). We say that  $\xi$ is $K$-bounded if  $|\xi| \le K$ with probability 1.

\begin{lemma} [Projection lemma,  \cite{TVdet} ] \label{TV1} 
Let $X= (\xi_1, \dots, \xi_n)$ be a random vector in $\C^n$  whose coordinates $\xi_i$ are independent $K$-bounded random variables with mean 0 and variance 1, where  $K \ge 10 (\E |\xi_i|^4 +1)$ for all $i$. Let $H$ be a subspace of dimension $d$ and $\Pi_H X$ be the length of the projection 
of $X$ onto $H$. Then
$$\P (| \Pi_H X - \sqrt{d}| \ge t ) \le 10 \exp(- \frac{t^2}{20 K^2}) . $$ 
\end{lemma} 

The projection lemma follows from the Talagrand's inequality (\cite[Chapter 4]{Ledoux}). The constants $10$ and $20$ are rather arbitrary. We make no attempt to optimize the constants in this paper. 

{\bf Notation.} We use standard  asymptotic notations such as $O, o, \Theta, \omega, \ll$ etc., under the assumption that $n \rightarrow \infty$. 
For a vector $X$, $\|X\|$ is its Euclidean norm and 
$\|X\|_{\infty}$ its infinity norm.  For a matrix $A \in \C^{n\times n}$,  $\| A\|_F $ and $\|A\|_2$ denote its Frobenius and spectral norm, respectively. 
All eigenvectors have unit length.

\subsection{Weighted projection lemmas}

Let us fix an orthonormal basis $\{u_1, \dots, u_d\}$ of $H$. We can express  $\Pi_H X $ as 
\begin{equation} \label{P1}  \Pi_H X =   (\sum_{i=1}^d |u_i^*  X|^2 )^{1/2}. \end{equation} 

In recent studies, we came up with  situations when the roles of the axes are not 
compatible. Formally speaking, one is required 
 to consider a weighted version of  \eqref{P1} where 
$ (\sum_{i=1}^d |u_i^*  X|^2 )^{1/2}$ is replaced by  $   (\sum_{i=1}^d c_i  |u_i ^*  X|^2 )^{1/2} $ with $c_i$ being non-negative numbers (weights).

We are able to prove a variant of  Lemma \ref{TV1} for this more general problem, which turns out to be useful in a number of applications, some of which will be discussed in the paper. 
Furthermore,  we can also weaken the assumption on  random vector $X$ in various ways. 

We say a random vector $X =(\xi_1, \dots, \xi_n) $ is $K$-{\it concentrated }(where $K$ may depend on $n$)  if there are 
constants $C, C' >0$ such that for any convex, $1$-Lipschitz function $F:  \C^n \rightarrow \R$ and any $t>0$
\begin{equation} \label{Tal} 
 \P ( |F(X) - M(F(X))| \ge t) \le C \exp (-C'  \frac{t^2}{K^2} ), 
\end{equation}  where $M(Y)$ denotes the median of a random variable $Y$ (choose an arbitrary one as there are many). 

Notice that the notion of $K$-concentrated is somewhat similar to the notion of threshold in random graph theory in the sense that 
if $X$ is $K$-concentrated then it is $c K$-concentrated for any constant  $c >0$. (Similarly, if $p(n)$ is a threshold for a  property $\CP$ (say, containing a triangle) then $cp(n)$ is also 
a threshold.) One can also replace the median by the expectation (see Lemma \ref{lemma:bound1}).  The dependence on $K$ on the RHS of \eqref{Tal} is flexible (one can replace $K^2$ by any function $f(K))$; however, the quality of the concentration bound will depend on $f(K)$
and we leave it as an exercise for the reader to work out this dependence.

{\it Examples of $K$-concentrated random variables} 
\begin{itemize}

\item If the coordinates of $X$ are iid standard gaussian (real or complex), then $X$ is  $1$-concentrated (see \cite{Ledoux}).

\item If $\xi_i$ are independent and $\xi_i$ are $K$-bounded for all $i$, then $X$ is $K$-concentrated (this is a corollary of Talagrand's inequality;
see \cite[Chapter 4]{Ledoux} or \cite[Theorem F.5]{TVsmall}). 

\item If $X$ satisfies the log-Sobolev inequality with parameter $K^2$, then it is $K$-concentrated (see \cite[Theorem 5.3]{Ledoux}). 

\item The coordinates $\xi_i$ of $X$ come from a random walk satisfying certain mixing properties (see \cite[Corollary 4]{Sam}; in this corollary $\| \Gamma \|$ plays the role of $K$). In this and the previous example, the coordinates of 
$X$ are not necessarily indepedent. 

\end{itemize}

\begin{lemma} \label{lemma:VW1} 
Let $X= (\xi_1, \dots, \xi_n)$ be a   $K$-concentrated 
 random vector in $\C^n$.
   Then there are constants $C, C' >0$ (which depend on, but could  be different from the constants 
in \eqref{Tal}) such that the following holds. Let $H$ be a subspace of dimension $d$ with an orthonormal basis $\{u_1,\ldots, u_d\}$. Then for any $1 \ge c_1, \dots, c_d \ge 0$
\begin{equation*}
\P \left( | \sqrt{\sum_{j=1}^d c_j |u_j^*X|^2} -\sqrt{\sum_{j=1}^d c_j}  | \ge t \right) \le C \exp(-C' \frac{t^2}{K^2} ).
\end{equation*}
\end{lemma}

In particular, by squaring, it follows that
\begin{equation} \label{ineq:VW} 
\P \left( | \sum_{j=1}^d c_j |u_j^*X|^2 -\sum_{j=1}^d c_j  | \ge 2t\sqrt{\sum_{j=1}^d c_j} + t^2 \right) \le C \exp(-C' \frac{t^2}{K^2} ).
\end{equation}

Another way to weaken the $K$-bounded assumption is to consider truncation. If $\xi$ is not bounded, but has light tail, then 
by setting $K$ appropriately, we can show that $\P( |\xi| \ge K)$ is negligible with respect to the probability bound we want to prove. 

Assume that the $\xi_i$ are independent with mean zero and variance one. Choose a number $K > 1$ and let 
$\eps_1 :=\max_{1 \le i \le n} \P( |\xi_i| >K)$. Set $\xi_i': =\xi_i \I _{|\xi_i| \le K}$  and let $\mu_i$ and $\sigma_i^2$ denote its mean and variance. Set $\eps_2 := \max_{1 \le i \le n} |\mu_i|$ and 
$\eps_3:= \max_{1\le i \le n} |\sigma_i^2-1 |$. Assume all $\eps_j \le 1/2$ (in practice this assumption is satisfied easily).

\begin{lemma} \label{lemma:VW2}  There are constants $C, C'>0$ such that the following holds. 
Let $X= (\xi_1, \dots, \xi_n)$ be a random vector in $\C^n$  whose coordinates $\xi_i$ are independent random variables with mean 0 and variance 1. 
Under the above notations, we have, for any $1 \ge c_1, \dots, c_n \ge 0$ and any $t >0$
\begin{equation} \label{VW2-0}  
\P \left( | \sqrt{\sum_{j=1}^d c_j |u_j^*X|^2} -\sqrt{\sum_{j=1}^d c_j}  | \ge t  + 4n^2 K^2 (\eps_2 +\eps_3) \right) \le C\exp(-C' \frac{t^2}{K^2} ) + n \eps_1 .
\end{equation}
\end{lemma}

\subsection{Concentration of random  quadratic forms}\label{section:Quad}

Consider a quadratic form $Y:= X^{*} A X$  where $X = (\xi_1, \dots, \xi_n)$ is, as usual,  a random vector 
 and $A= (a_{ij}) _{1 \le i, j \le n}$  a deterministic matrix. As application of the  new projection lemmas,  we prove  a large deviation result for  $Y$, which can  be seen as the quadratic version of the standard Chernoff bound.


\begin{theorem}[Concentration of quadratic forms  I] \label{thm:quadratic}  Let $X$ be a  $K$-concentrated random vector in $\C^n$. 
Then there are constants $C, C' >0$ such that  for any matrix $A$
 \begin{equation}\label{eq:wp-hw1}
\P(|X^*A X - \trace(A)| \ge t) \le C  \log n \exp(- C' K^{-2}  \min\{\frac{t^2}{\|A\|_F^2 \log n }, \frac{t}{ \| A\| _2}  \}).
\end{equation}
\end{theorem}

 \begin{theorem} [Concentration of quadratic forms II]  \label{thm:quadratic2}  Let $X$  and $\eps_1, \eps_2, \eps_3$ be as in Lemma \ref{lemma:VW2}. There are constants $C, C' >0$ such that the following holds. 
 Assume  $n^2 K^2   (\eps_2 +\eps_3)  =o(1)$, then 
$$  \P(|X^*A X - \trace(A)| \ge t ) \le   C  \log n \exp(- C' K^{-2}  \min\{\frac{t^2}{\|A\|_F^2 \log n }, \frac{t}{ \| A\| _2}  \}) + n \eps_1. $$
\end{theorem}

As an illustration, let us consider 
 the case when $\xi_i$ are sub-exponential.  We say that $\xi$ is sub-exponential with exponent $\alpha >0$ (with accompanying positive constants $a$ and $b$) if 
 for any $ t >0$ 
$$ \P( |\xi -\E \xi| \ge t^{\alpha} ) \le a \exp(-bt ) . $$

 \begin{corollary} [Concentration of quadratic forms with sub-exponential variables]  \label{cor:quad} 
 Assume that  $\xi_i$ are independent and sub-exponential (with exponent $\alpha >0$) random variables with mean 0 and variance 1. Then  
 there are constants $C, C' >0$ such that for any $t$ satisfying 
 \begin{equation}  \label{t} 
 t =\omega ( ( \| A\|_F + \log^{\alpha} n \|A\|_2 ) \log^{\alpha+1} n), 
 \end{equation}
 we have 
  \begin{equation} \label{VW5}  \P ( |X^*AX - \trace(A)|  \ge t )  \le C \exp(-C'   \min \{ (\frac{t}{ \|A\|_F \sqrt {\log n}} )^{\frac{1}{ \alpha +1/2 }}, (\frac{t}{\|A\|_2 }) ^{\frac{1}{ 2 \alpha + 1}} \} )  .\end{equation}

 \end{corollary}

  Quadratic forms of random variables appear frequently in probability  and  the large deviation problem has been considered by several researchers, with 
   first and perhaps most famous result being  by  Hanson-Wright inequality  \cite{HW}.  In many cases, our results improve earlier  results significantly; see Section \ref{section:Qproof} for more details.

\subsection{Norm of random eigenvectors} \label{section:norm}

Let $M_n$ be a symmetric $\pm 1$ matrix (the upper diagonal and diagonal entries are iid  Bernoulli random variables taking values $\pm 1$ with probability $1/2$). This is an important object in both probabilistic combinatorics and the theory of random matrices.  Let $u $ be an arbitrary unit eigenvector of $M_n$. We investigate  the following 
natural question, 
\vskip2mm

\centerline {\it How big is $\| u\|_{\infty}$  ? } 

A good bound on the infinity norm of the eigenvectors is
 important   in spectral analysis of graphs  and many other applications, such as the studies of nodal domains  (see for instance \cite{DLL} and the references therein). 
 Recently, it  plays a crucial role in breakthrough works concerning   local statistics of random matrices  (see Section \ref{section:LSCL} and also  \cite{Erdossur, TVsur} for surverys). 
 
 Set $W_n  =\frac{1}{\sqrt n} M_n$.  Thanks to the classical Wigner's semi-circle law (see Section \ref{section:LSCL}), we know that most of the eigenvalues of $W_n$ belong to the interval $(-2 +\epsilon, 2 -\epsilon)$. 
 Using our new concentration results, we are able to obtain (what we believe to be) the optimal bound for the eigenvectors corresponding to these eigenvalues.

  \begin{theorem}[Infinity norm of eigenvectors]  \label{theorem:delo} Let $M_n$ be a $n\times n$ symmetric matrix whose upper diagonal entries are iid random variables that takes values $\pm 1$ with the same probability. Let $W_n  =\frac{1}{\sqrt n} M_n$. For any constant $C_1 >0$, there is a constant $C_2 >0$ such that the following holds.

\begin{itemize}
\item (Bulk case) With probability at least $1- n^{-C_1}$, for  any fixed  $\epsilon >0$ and any $1\le i \le n$ with $\lambda_i(W_n) \in [-2+\epsilon, 2-\epsilon]$, there is a unit eigenvector  $u_i(W_n)$ of $\lambda_i (W_n)$ satisfying 
 $$\| u_i(W_n) \|_{\infty} \le \frac{C_2  \log^{1/2} n}{\sqrt n}.$$

\item (Edge case) With probability at least $1- n^{-C_1}$, for any $\epsilon >0$ and any $1\le i \le n$ with $\lambda_i(W_n) \in [-2-\epsilon, -2+\epsilon] \cup [ 2-\epsilon, 2+\epsilon]$,  there is a unit eigenvector  $u_i(W_n)$ of $\lambda_i (W_n)$ satisfying 
 $$\| u_i(W_n) \|_{\infty} \le \frac{C_2  \log  n}{\sqrt n}.$$
\end{itemize}
\end{theorem}

 The best previous bound was of the form $\frac{\log^ C n}{ n^{1/2} }$ for some large (usually not explicit) constant $C$ \cite{ESY2, ESY1, ESYweg, TVsmall}.
 We conjecture that the bound $O( \sqrt {\log n/n})$ is sharp (it is easy to see that this is the case if the entries are standard gaussian) and also that it  holds for all eigenvectors. Very recently, Rudelson and Vershynin \cite{RVvector} also studied the norm of random eigenvectors using a geometric method,  which is different from our approach discussed in Section \ref{section:LSCL}.

 \subsection{The local semi-circle law} \label{section:LSCL}

 Denote by $\rho_{sc}$ the semi-circle density function with
support on $[-2,2]$,
\begin{equation}\label{semi}
 \rho_{sc} (x):= \begin{cases} \frac{1}{2\pi} \sqrt {4-x^2}, &|x| \le 2 \\ 0,
&|x| > 2. \end{cases} 
\end{equation}

Let us recall the classical Wigner's semi-circle law:

\begin{theorem}[Semi-circular law]\label{theorem:Wigner} Let  $M_n$ be a random Hermitian matrix whose entries on and above the diagonal are iid bounded random variables with zero mean and unit variance  and $W_n= \frac{1}{\sqrt{n}} M_n$.  Then for any real number $x$,
$$\lim_{n \rightarrow \infty} \frac{1}{n} |\{1 \leq i \leq n: \lambda_i(W_n) \le  x \}|
 = \int_{-2}^x \rho_{sc}(y)\ dy$$
in the sense of probability, where we use $|I|$ to denote the cardinality of a finite set $I$.  
\end{theorem}

The key tool for bounding the infinity norm of  an  eigenvector  is a statement of the following type: any  interval of length at least $T$ (which tends  to zero with $n$) 
 in the spectrum $[-2,2]$ contains  an eigenvalue, with high probability. The 
quality of the bound will depend on how small $T$ is.  This approach was developed by Erd\H os, Schlein and Yau in  \cite{ESY2, ESY1, ESYweg}, leading to eigenvector norm bounds 
of  order $n^{-2/3}, n^{-3/4}$ and finally $n^{-1 +o(1)}$. A simpler argument, following the same approach, was developed by Tao and the first author in \cite{TVsmall} (see \cite[Section 4]{TVsmall} for  a problem 
concerning random non-hermitian matrices).

One  way to attack the above problem is 
to show that the semi-circle law holds for small intervals (or at small scale). Intuitively, we would like to have with high probability that
$$| N_I - n \int_I  \rho_{sc}(x) \ dx | \le \delta n |I|  , $$ for any interval $I$ and fixed $\delta >0$, where $N_I$ denotes the number of eigenvalues of $W_n:= \frac{1}{\sqrt n}M_n$ in the interval $I$.  Of course,  the 
reader can easily see that $I$ cannot be 
arbitrarily  short (since $N_I$ is an integer).  Following \cite{ESYweg}, we call a  statement of this kind a local semi-circle law (LSCL).

A natural question arises:  {\it how short can $I$ be ? } 
Formally, we say that the LSCL  holds at a scale $f(n)$ if   
with probability $1-o(1)$ $$| N_I - n \int_I  \rho_{sc}(x) \ dx | \le \delta n |I|  , $$ for any interval $I$ in the bulk  of length $\omega( f(n))$ and any fixed $\delta >0$. Furthermore, we say that 
$f(n)$ is a {\it threshold scale} if the LSCL holds at scale $f(n)$ but does not holds at scale $g(n)$ for any function $g(n)=o(f(n))$. (The reader may notice a
 similarity between this definition and the definition of threshold functions for random graphs.)   We would like to raise the following problem. 

\begin{problem} \label{question:main0} 
Determine the threshold scale (if exists). 
\end{problem} 

We do not know a  sharp estimate for the threshold for any matrix ensembles,  even in the basic GUE  (random matrix with complex gaussian entries) and GOE  (random matrix with real gaussian entries) cases. A recent result  by Ben Arous and Bourgade \cite{Ben} shows that the maximum gap between two  consecutive 
(bulk) eigenvalues  of GUE  is of order $\Theta (\sqrt{\log n}/n)$, with high probability. Thus, if we partition the bulk into intervals of length $\alpha \sqrt{\log n}/n$ for a
sufficiently  small $\alpha$, one of these intervals contains at most one eigenvalue. Therefore,   we   expect that in natural ensembles, the LSCL does not  hold below the $\sqrt{\log n}/n$ scale.
In \cite{ESYweg, TVuniv}, upper bound of the form $\log^C n/n$ was proved for some large value of $C$. Here we are going to show

\begin{theorem}[Threshold for local semi-circle law]  \label{theorem:main} 
 Let $M_n$ be a Hermitian  matrix whose upper diagonal entries are independent random variables  with mean $0$ and variance $1$. Assume furthermore that  for $1\le i \le n$, the vectors $X_i$, obtained by deleting the $i$-th entry of the $i$-th row vector of $M$, are $K$-concentrated. Then the threshold scale for LSCL is bounded from above by 
$K^2 \log n /n$. 
\end{theorem}

In the GUE case, the gap between the upper and lower bound is only $O(\sqrt {\log n})$ and it is an  intriguing problem to remove this factor. 
We also conjecture that  Ben Arous and Bourgade's result  on the largest gap holds for  $\pm 1$ random matrices. 

The results of Section \ref{section:norm} and Section \ref{section:LSCL} also hold for random sample covariance matrices. We sketch the results and proofs in the appendices. 

\vskip2mm

{\it Structure of the paper.} In the next section, we prove the new projection lemmas. In Section \ref{section:Qproof}, we prove the new concentration inequalities for quadratic forms and make a comparison with prior results. 
The next section, Section \ref{section:fact} can be seen as a preparation step in which we recall facts about random matrices. We prove the new threshold for the local semi-circle law in Section \ref{section:LSCLproof}, and the bound 
on the infinity norm of eigenvectors in Section  \ref{section:evproof}. The appendices contain proofs concerning random sample covariance matrices. 

{\bf Acknowledgement.}
The authors would like to thank the anonymous referees for their careful reading and constructive suggestions.

 \section{Proofs of Lemmas \ref{lemma:VW1} and \ref{lemma:VW2} } 
 
 {\it Proof of Lemma \ref{lemma:VW1}.} 
Set  $f(X) :=\sqrt{\sum_{j=1}^d c_j |u_j^*X|^2}$. Thus, $f$ is a function from $\C^n$ to $\R$.  

We first observe that $f(X)$ is convex. Indeed, for $0\le \lambda, \mu \le 1$ with $\lambda+\mu=1$ and any $X,Y \in \BBC^n$, by Cauchy-Schwardz inequality,
\begin{equation*}
\begin{split}
f(\lambda X+ \mu Y) &\le \sqrt{\sum_{j=1}^d c_j (\lambda|u_j^*X|+\mu|u_j^*Y|)^2}\\
&\le \lambda\sqrt{\sum_{j=1}^d c_j |u_j^*X|^2}+\mu\sqrt{\sum_{j=1}^d c_j |u_j^*Y|^2}=\lambda f(X)+\mu f(Y).
\end{split} 
\end{equation*}

Next, we show that  $f(X)$ is {1-Lipschitz}. Notice that $f(X) \le  \sqrt{\sum_{j=1}^d |u_j^*X|^2} \le \|X\|$. Since $f(X)$ is convex, one has $$\frac{1}{2}f(X)=f(\frac{1}{2}X)=f(\frac{1}{2}(X-Y)+\frac{1}{2}Y)\le \frac{1}{2}f(X-Y)+\frac{1}{2}f(Y).$$ 

Thus $f(X)-f(Y)\le f(X-Y)$ and $f(Y)-f(X) \le f(Y-X)=f(X-Y)$, which imply 
$$|f(X)-f(Y)|\le f(X-Y) \le \|X-Y\|.$$

Thus, by the definition of $K$-concentrated property,
\begin{equation}\label{eq:median}
\Pr(|f(X)-M(f(X))| \ge t) \le C \exp(- C' \frac{t^2}{K^2} )
\end{equation} for some constants $C, C' >0$. 

To conclude the proof, it suffices to show 
$|M(f(X))-\sqrt{\sum_{j=1}^d c_j}| =O(K).$ We use the following lemma. The proof of the lemma is classical and thus omitted.

\begin{lemma}  \label{lemma:bound1} Let $Y$ be a real random variable. 
Assume $\P( |Y- \mu| \ge t ) \le f(t)$, where $\int_0 ^{\infty} f(x) dx = O(1)$, 
then $|\E Y - \mu| = O(1)$.  Assume furthermore that  $Y$ is non-negative, $\mu \ge 0$, $ \sigma = \sqrt{\E Y^2}$,  and  $\int_0^{\infty} x f(x) dx= O(1)$, then $| \E Y -\sigma| = O(1) . $

\end{lemma} 

To apply this lemma, set $c_i'  := \frac{c_i}{ \max _{1 \le i \le d } c_i }$, 
 $Y := \frac{1}{K} \sqrt {\sum_{i=1}^d c'_i |u_i^* X|^2 } $ and $\mu:= M(Y)$.
 We have, by the $K$-concentration property
 $$\P (|Y -\mu|  \ge t) =  \P( | \sqrt {\sum_{i=1}^d c'_i |u_i^* X|^2 }-  M (  \sqrt {\sum_{i=1}^d c'_i |u_i^* X|^2 } ) | \ge tK ) \le C \exp(- C' t^2). $$
 
 Set  $f(x)= C \exp( -C' x^2)$.  The assumptions on $f(x)$ in Lemma \ref{lemma:bound1}  are trivially satisfied. 
As $X$ is isotropic, $\sigma^2 =\E Y^2 = \frac{1}{K^2} \sum_{i=1}^d c_i'$. 
It follows from Lemma \ref{lemma:bound1} that 
$$M(Y) =   \frac{1}{K} \sqrt {\sum_{i=1}^d c_i' } +O(1). $$ 
 
Renormalizing, we obtain
$$M( \sqrt {\sum_{i=1}^d c_i |u_i^* X|^2}) =   \sqrt{\sum_{i=1}^d c_i} + O(K \sqrt{\max_{1 \le i \le d} c_i} ), $$  
which concludes the  proof of Lemma \ref{lemma:VW1}.

\vskip2mm 

{\it Proof of Lemma \ref{lemma:VW2}.}  Recall that $\xi_i'=\xi_i \I _{|\xi_i| \le K}$ has mean $\mu_i$ and variance $\sigma_i^2$. The parameters $\eps_1 =\max_{1 \le i \le n} \P( |\xi_i| >K)$, $\eps_2 = \max_{1 \le i \le n} |\mu_i|$ and  $\eps_3= \max_{1\le i \le n} |\sigma_i^2-1 |$ satisfy all $\eps_j \le 1/2$. Define $\tilde \xi_i:= \frac{\xi_i' -\mu_i}{\sigma_i}$. The $\tilde \xi_i$ are independent with mean zero and variance 1 and are $2K$-bounded.  Let $X':= (\xi_1', \dots, \xi_n')$ and $\tilde X:= (\tilde \xi_1, \dots, \tilde \xi_n)$. 
 It is obvious that 
 $$
 \P \left( | \sqrt{\sum_{j=1}^d c_j |u_j^*X|^2} -\sqrt{\sum_{j=1}^d c_j}  | \ge t \right)  \le \P \left( | \sqrt{\sum_{j=1}^d c_j |u_j^*X'|^2} -\sqrt{\sum_{j=1}^d c_j}  | \ge t \right)  + n\eps_1. 
 $$
 
 The next observation   is that if $\eps_2, \eps_3$ are small, then $\sum_{1 \le i \le d} c_i |u_i^* X'|^2$ and   $\sum_{1 \le i \le d} c_i |u_i^* \tilde X|^2$ are more or less the same. By definition, we have 
 with probability one
 $$
 |\xi_i' -\tilde \xi_i | = | \frac{\xi_i'(\sigma_i-1) + \mu_i }{\sigma_i} | \le 2 (K \eps_3+ \eps_2). 
 $$
 
 It follows that  $D := X'  -\tilde X$ has norm at most $ 2 n^{1/2} (K \eps_3+ \eps_2)$ with probability one. On the other hand, 
$$ 
\Big| \sum_{1 \le i \le d} c_i |u_i^* X'|^2 - \sum_{1 \le i \le d} c_i |u_i^* \tilde X|^2 \Big| \le 2 \sum_{1\le i\le d} c_i  |u_i^*X_i'|  |u_i^* D| +  \sum_{1\le i \le d}c_i | u_i^* D |^2. 
$$ 

As $u_i$ are unit vectors,  $|u_i^* X_i'| \le \| X_i'\| \le \sqrt n K$ and $| u_i^*D_i| \le \| D_i\| \le 2 \sqrt n (K \epsilon_2+ \epsilon_3) $ (these bounds are 
 generous and can be improved by a polynomial factor in certain cases, but in applications such improvement rarely matters). It follows, again rather generously, 
 $$ |\sum_{1 \le i \le d} c_i |u_i^* X'|^2 - \sum_{1 \le i \le d} c_i |u_i^* \tilde X|^2| \le 4n  \sum_{i=1}^d c_i K^2 (\eps_2 +\eps_3)
  \le 4n^2K^2 (\eps_2 +\eps_3).  
  $$
  
Applying Lemma \ref{lemma:VW1} for $\tilde X$, we obtain Lemma \ref{lemma:VW2}.

\vskip2mm 
In practice, $\eps_j$ are typically super-polynomially small, i.e. $n^{-\omega(1)}$,  which yields  $4n^2 K^2 (\eps_2 +\eps_3) =o(1)$. This term can be ignored
(by slightly changing the values of $C, C'$ if necessary) and we end up with a more friendly inequality 
\begin{equation} \label{VW2-1}  
\P ( | \sqrt{\sum_{j=1}^d c_j |u_j^*X|^2} -\sqrt{\sum_{j=1}^d c_j}  | \ge t ) \le C\exp(-C' \frac{t^2}{K^2} ) + n \eps_1. 
\end{equation}

In the sub-exponential case, 
for a sufficiently large $K$ (compared to $a$ and $b$), 
$\eps_j \le  \exp(- \frac{b}{2} K^{1/\alpha})$ for $j=1,2,3$.  For $K= \omega (\log^{\alpha} n)$, 
 $ n^2K^2\exp( -\frac{b}{2}  K^{1/\alpha}) =o(1)$  and 
\eqref{VW2-1} yield
\begin{equation} \label{VW2-2}  
\P \left( | \sqrt{\sum_{j=1}^d c_j |u_j^*X|^2} -\sqrt{\sum_{j=1}^d c_j}  | \ge t   \right) \le C \exp(- C' \frac{t^2}{K^2} ) + n\exp(- \frac{b}{2} K^{1/\alpha})  .
\end{equation}

 \section {Random Quadratic Forms } \label{section:Qproof}

\subsection{Proofs of new results} 

Let us first prove Theorem \ref{thm:quadratic}. 
 Notice that if $Y = X^* A X$, then $Y +\bar Y = X^*(A+A^*) X$ and 
$Y- \bar Y = X^* (A- A^* ) X$.  Since 
$$
Y - \trace A = \frac{1}{2} [ (Y+\bar Y) - \trace (A +A^*) ]  +  \frac{1}{2}  [ (Y - \bar Y ) -\trace (A-A^*) ] , 
$$  
we have 
$$
\P( | Y -\trace A | \ge t) \le \P (|(Y+\bar Y) - \trace (A +A^*)| \ge t )  + \P ( |\sqrt{-1}(Y - \bar Y ) -\trace (\sqrt{-1}(A-A^*) ) | \ge t ). 
$$

Moreover,  as $\| A + A^{*}\|_F, \| A-A^*\| _F =O( \|A\|_F)$ and $\| A +A^*\| _2, \|A-A^*\| _2 = O(\| A\|_2)$, it suffices to prove the theorem in the case $A$ is Hermitian. 

Next, we observe that any  Hermitian matrix $A$ can be written as $A:= A_1 -A_2$ where $A_i$ are positive semi-definite and 
$\max_{i=1,2} \|A_i\|_2 \le \|A\|_2, \max_{i=1,2} \|A_i\|_F \le \|A\|_F$. (In fact, the positive eigenvalues of $A_1$ are the positive eigenvalues of  $A$ and the positive eigenvalues of 
$A_2$ are the absolute values of the negative eigenvalues of $A$.) 
This enables us to further reduce the problem to  the case when $A$ is positive semi-definite.

Finally,  as the content of the theorem is invariant under scaling,  we can assume that $\|A\|_2=1$. Let $1= c_1  \ge  c_2, \dots, c_n \ge 0$ be the eigenvalues of $A$
together with corresponding orthonormal eigenvectors  $\{u_1,\ldots,u_n\}$. We have 
\begin{equation}  \label{general}   X^*AX- \tr(A) = \sum_{j=1}^n c_j |u_j^*x|^2 - \sum_{j=1}^n c_j. \end{equation} 

This is precisely the setting of the projection lemmas.  By Lemma \ref{lemma:VW1}, we know that 
for any numbers $0 \le d_j  \le 1$, $j \in J$, 
\begin{equation}\label{eq:WP1}
\P(|\sum_{j\in J} d_j |u_j^*X|^2 - \sum_{j\in J} d_j | \ge 2t \sqrt{ \sum_i d_i} +t^2) \le C \exp(-C' K^{-2} t^2 ). 
\end{equation}

However, it is somewhat  wasteful to apply this 
directly to \eqref{general}. We will perform  an extra partition step. Set 
$$J_k := \{1\le j\le n : \frac{1}{4^{k+1} } \le c_j \le \frac{1}{4^k } \}, 0 \le k \le k_0:= 10 \log n, $$ and  let $J_{k_0+1}$ be the collection of the remaining indices.

For each $0\le k \le k_0+1$, apply Lemma \ref{lemma:VW1} to  $d_i := 4^{k} c_i, c_i \in J_k$, we have, for any $s \ge 0$
$$
\P ( | \sum_{i \in J_k}  4^k  c_i( | u_i^* X|^2 -1)|  \ge 2s  \sqrt{\sum_{i \in J_k} 4^k c_i } + s^2 ) \le C \exp(- C' K^{-2} s^2). 
$$
 
Set $s : = \frac{t}{\|A\|_F}$ and simplify by $4^k$,  the above inequality becomes 
$$
\P( | \sum_{i \in J_k} c_i( | u_i^*X|^2 -1)| \ge \frac{2t}{2^k \|A\|_F } \sqrt {\sum_{i \in J_k} c_i }  + \frac{t^2}{4^k \|A\|_F ^2 })  \le C \exp(- C' K^{-2} \frac{t^2}{ \|A\|_F^2 } ) . 
$$

Apparently,  $\sum_{k=0}^{k_0+1} \frac{t^2}{4^k \|A\|_F ^2 } \le 2 \frac{t^2}{ \|A\|_F^2 } $. Moreover, 
$\sum_{i \in J_{k_0+1}}  c_i \le n \times n^{-5} = n^{-4} $ and 
\begin{eqnarray*}  
\sum_{ 0 \le k \le k_0} 2^{-k} \sqrt {\sum_{i \in J_k} c_i } &\le&  k_0^{1/2}  (\sum_{k=0}^{k_0} 4^{-k} \sum_{i \in J_k} c_i )^{1/2} \\
&\le&8 \log^{1/2} n (\sum_{k=0}^{k_0} \sum_{i \in J_k} c_i^2 )^{1/2}  \\
&\le& 8 \log^{1/2} n  \| A\|_F, 
\end{eqnarray*}  
by Cauchy-Schwartz inequality.

Putting the above estimates together and using the union bound, we obtain 
$$ 
\P (| \sum_{i=1}^n c_i (| u_i^* X|^2-1) | \ge 16  \log^{1/2} n  t + 2 \frac{t^2}{\|A\|_F^2 } + n^{-2} ) \le  C \log n  \exp(- C'K^{-2} \frac{t^2}{\|A\|_F^2 }). 
$$

We  can ignore the small term $n^{-2}$ (by slightly adjusting the constant 16),
the desired bound follows. 

\begin{remark} If we have more information about $A$, the 
 $\log n$ term can be improved. For instance, if all eigenvalues of $A$ are comparable, then 
we do not need  this term.
 \end{remark}
The proof of  Theorem \ref{thm:quadratic2} uses Lemma \ref{lemma:VW2} and is left as an exercise. 
To prove Corollary \ref{cor:quad}, notice that we can  obtain an 
  analogue 
 of \eqref{VW2-2} 
 \begin{equation} \label{VW4}  \P ( |X^*AX - \trace(A)| \ge t  ) \le  C \exp(- C' K^{-2}  \min\{\frac{t^2}{\|A\|_F^2 \log n }, \frac{t}{ \| A\| _2}  \})  + n\exp(- \frac{b}{2} K^{1/\alpha}) , \end{equation}
\noindent under the assumption that $K =\omega ( \log^{\alpha} n)$. 

 To optimize the bound, we choose  $K$ such that 
  $K^{-2}  \min\{\frac{t^2}{\|A\|_F^2 \log n }, \frac{t}{ \| A\| _2}  \} = K^{1/\alpha}$. This leads to 
 setting $K:= \min \{ (\frac{t}{ \|A\| _F  \sqrt {\log n}} )^{\frac{2}{2 +1/\alpha}}, (\frac{t}{\|A\|_2 }) ^{\frac{1}{ 2+ 1/\alpha}}  \}$. 
 Assume 
 \begin{equation}  
 \label{t1} t =\omega ( ( \| A\|_F + \log^{\alpha} n \|A\|_2 ) \log^{\alpha+1} n). 
 \end{equation}
 
 This assumption guarantees  $K =\omega ( \log^{\alpha} n)$. It also implies 
 $n\exp( -\frac{b}{2}  K^{1/\alpha}) \le \exp( -\frac{b}{3} K^{1/\alpha})$, proving Corollary \ref{cor:quad}.

\subsection{Comparison to earlier results}

  In 1971, Hanson and Wright \cite{HW} obtained the first  important inequality   for sub-gaussian random variables.
  
\begin{theorem}[Hanson-Wright inequality] \label{HW}
 Let  $X=(\xi_1,\ldots, \xi_n) \in \R^n$ be  a random vector with $\xi_i$ 
 being  iid  symmetric  and sub-gaussian random variables with mean 0 and variance 1. 
 There exist constants $C,C'>0$ 
  such that the following holds. Let $A$ be a real matrix of size $n$ with entries $a_{ij}$ and 
 $B :=(|a_{ij}|)$. Then 
 \begin{equation}\label{eq:HW}
\P(|X^T A X -\tr(A)| \ge t) \le C\exp (-C'\min\{\frac{t^2}{\|A\|_F^2},  \frac{t}{\|B\|_2 } \} )
\end{equation}
for any $t>0$.
\end{theorem}

Later, Wright \cite{Wright} extended Theorem \ref{HW} to non-symmetric random variables.  Recently, Hsu, Kakade and Zhang \cite{HKZ} showed that one can obtain 
a better upper tail (notice that $\|B\|_2$ is replaced by $\|A\|_2$) 
\begin{equation}\label{eq:HKZ}
\P(X^T A X -\tr(A) \ge t) \le C\exp (-C'\min\{\frac{t^2}{\|A\|_F^2},  \frac{t}{\|A\|_2} \} )
\end{equation}
\noindent under a considerably  weaker assumption (which, in particular, does not require 
the $\xi_i$ to be independent). On the other hand, their method does not cover the lower tail. 
 Let us pause here  to point out  a strong distinction from the linear case and the quadratic case: In the linear case (Chernoff  type bounds), the lower tail follows from the upper tail 
 by simply switching $\xi_i$  to $-\xi_i$, but this trick is useless in the quadratic case. Recently, Rudelson and Vershynin \cite{RVhw} proved the Hanson-Wright inequality 
 \begin{equation}
\P(|X^T A X -\tr(A)| \ge t ) \le C\exp (-C'\min\{\frac{t^2}{\|A\|_F^2},  \frac{t}{\|A\|_2} \} ),
\end{equation}
assuming $\xi_i$ are sub-gaussian.

In the previous papers, the random variables $\xi_i$ are required to be real. 
Few years ago,   motivated by the delocalization problem for random matrices, 
Erd\H os, Schlein and Yau  \cite{ESYweg}   considered the complex case.  
By assuming either both the real and imaginary parts of $\xi_i$ are iid sub-gaussian or the distribution of $\xi_i$ is rotationally symmetric (real and imaginary parts still sub-gaussian),  they proved 
\begin{equation}  \label{ESY1}   \P(|X^* A X -\tr(A)| \ge t) \le C\exp( -C' \frac{t}{\|A\|_F} ). 
\end{equation} 

Later,   
Erd\H os, Yau and Yin \cite{EYYbulk}  showed that if $\xi_i$ are independent sub-exponential random variables with exponent $\alpha >0$, having mean 0 and variance 1, then 
\begin{equation}   \label{ESY2}  \P(|X^* A X -\tr(A)| \ge t) \le C \exp  (-C'  
(\frac{t}{\|A\|_F } )^{\frac{1}{2 + 2 \alpha}  } ). 
\end{equation} 

To simplify the  comparison, let us  ignore the $\log n$ terms in our theorems (which play little role  in practice). 
If $K= O(1)$, then the main difference between Theorem \ref{HW} of Hanson and Wright 
 and Theorem \ref{thm:quadratic} is that the term $\| B\|_2$ in  Theorem \ref{HW} is now replaced by $\|A\|_2$. It is easy to see that 
 $\| B\| _2 \ge \|A\|_2 $ for any real matrix $A$. In fact, in many cases, $\|B\|_2$ is significantly larger than $\|A\|_2$. For instance, 
 a random matrix  $A$ with entries of order 1 typically has spectral  norm of order $\sqrt n$, but in  this case it is clear that $\|B\|_2$ has
 spectral  norm of order 
 $n$ (as all row sums are of this order). The same holds for several classical explicit  matrices, such as the Hadamard matrix. In these cases, our bound improves
 Hanson-Wright's significantly. Furthermore, our result applies in the complex case while the approach used 
 by Hanson and Wright is restricted to the real case. 
 
 Comparing to \eqref{ESY1},  we do not need the fairly restricted assumption that   either both the real and imaginary parts of $\xi_i$ are iid sub-gaussian or the distribution of $\xi_i$ is rotationally symmetric.
 In the case $K=O(1)$, both terms $\frac{t^2}{\|A\|_F^2 \log n }$ and $ \frac{t}{ \| A\| _2}$ in our bound can be considerably 
  larger than $\frac{t}{ \|A\|_F}$.  For instance, $ \frac{t}{ \| A\| _2}$ 
 and $\frac{t}{\| A\|_F}$  differ by a factor  $\sqrt n$ in both the random and Hadamard cases.

In  order to make a Hanson-Wright type bound 
 non-trivial, we need to assume  $t \ge   \| A\|_F +\|A\|_2 $. In many applications,  we want the probability bound to be 
 polynomially or even super-polynomially small, i.e. $n^{-O(1)}$ or $n^{-\omega(1)}$. This  requires  a  lower bound 
 $ \log^{\Omega (1)} n  ( \| A\|_F +\|A\|_2 )$ on $t$,  which is consistent with the assumption \eqref{t} in Corollary \ref{cor:quad}.

 Notice that \eqref{VW5} compares favorably to \eqref{ESY2}. For the term $ \frac{t}{ \|A\| _F } $, the exponent $\frac{1}{\alpha +1/2}$ is superior to  $\frac{1}{2\alpha +2 }$
 (notice that  we are talking about  a double exponent, so an improvement here could improve the quality of the bound quite a lot). 
  For the term $\frac{t}{\|A\|_2 }$, the exponent $\frac{1}{2 \alpha +1} $ is still better than $\frac{1}{2 \alpha +2}$. Furthermore, $\| A\|_2$ can be 
significantly smaller than $\|A\| _F$, as discussed earlier.

\section{Random matrices and the Stieltjes transform} \label{section:fact} 

This section serves as a preparation, in which we recall several facts about  random matrices. 
The empirical spectral distribution (ESD) function of  the  $n\times n$ Hermitian matrix 
$W_n:= \frac{1}{\sqrt n} M_n = \frac{1}{\sqrt n} (\zeta_{ij})_{1\le i, j \le n} $ is a one-dimensional function $$F^{\bf W_n}(x)=\frac{1}{n} |\{ 1\le j \le n: \lambda_j(W) \le x\}|,$$ where $|\mathbf{I}|$ denotes the cardinality of a set $\mathbf{I}$.
We are going to focus on the case when the entries of $M_n$ are $K$-bounded; it is easy to extend this assumption to  $K$-concentrated.

The Stieltjes transform of a real measure $\mu(x)$ is defined for any complex number $z$ not in the support of $\mu$ as 
$$s(z) =\int _{\R}  \frac{1}{x-z} d \mu(x). $$

Thus, the   Stieltjes transform $s_n(z)$ of  $W_n$ is 
$$\displaystyle s_n(z) =\int_{\mathbb{R}} \frac{1}{x-z} dF^{\bf W_n}(x) = \frac{1}{n} \sum_{i=1}^{n} \frac{1}{\lambda_i(W_n)-z}.$$

Furthermore, the Stieltjes transform $s_{sc}(z)$ of the semi-circle distribution is   $$ s_{sc}(z):=\int_{\mathbb{R}} \frac{\rho_{sc} (x) }{x-z} dx= \frac{-z+\sqrt{z^2-4}}{2},$$ where $\sqrt{z^2-4}$ is the branch of square root with a branch cut in $[-2,2]$ and asymptotically equals $z$ at infinity \cite{BS}. 

The beauty (and power) of the Stieltjes transform lies in the fact  that it has   a clear linear algebra content; $s_n(z)$ of $W_n$ is exactly  the trace of the 
matrix $(W_n -zI) ^{-1} $. This allows us to compute the Stieltjes transform by looking at the diagonal entries of $(W_n -zI) ^{-1} $. In matrix theory, Stieltjes transform  
plays the role Fourier transform in analysis. If  the Stieltjes transforms of two spectral measures  are close to each other (for all $z$), then 
the two measures are more or less the same. In particular,   if $s_n(z)$ is close to $s_{sc}(z)$, then the spectral distribution of 
$W_n$ is close to the semi-circle distribution (see for instance \cite[Chapter 11]{BS}, \cite{ESY1}).  
 We are going to  use the following lemma. 

\begin{lemma} \label{lemma1} Let $M_n$ be a random Hermitian matrix with independent $K$-bounded entries with mean 0 and variance 1. 
 Let $1/n < \eta< 1/10$ and $L, \varepsilon, \delta >0$. For any constant $ C_1>0$, there exists a constant $C >0$ such that if one has the bound $$|s_n(z)-s_{sc}(z)| \le \delta$$ with probability at least $1-n^{-C}$ uniformly for all $z$ with $|\text{Re}(z)| \le L$ and $\text{Im}(z) \ge \eta$, then for any interval $I$ in $[-L+\varepsilon, L-\varepsilon]$ with $|I| \ge \text{max}(2\eta, \frac{\eta}{\delta} \log \frac{1}{\delta})$, one has $$|N_I- n \displaystyle\int_I {\rho}_{sc}(x)\,dx| \le \delta n |I|$$ with  probability at least $1-n^{-C_1}$. 
\end{lemma}

This  is  \cite[Lemma 64]{TVuniv}, which, in turn, is a variant of \cite[Corollary 4.3]{ESY1}. 

An appropriate application of Lemma \ref{lemma1} will imply Theorem \ref{theorem:main}. (As a matter of fact, we a going to prove a little bit more.) 
In order to use this lemma, we set  $L=4, \varepsilon=1$, and critically 
 $$\eta :=\frac{K^2 C^2\log n }{n\delta^6},$$ where $C = C_1  +10^4$.  We are going to show that  
\begin{equation}\label{eq:region}
|s_n(z)-{s}_{sc}(z)| =o( {\delta})
\end{equation}
holds with probability  at least $1-n^{-C}$  for any fixed  $z$ in the region $\{z \in \mathbb{C}: |\text{Re}(z)| \le 4$, $\text{Im}(z) \ge {\eta} \}$. Notice that in this statement we fix $z$. However, it is simple to 
strengthen the statement to hold for all $z$, using an  $\epsilon$-net argument,  exploiting the fact that  $s_n(z)$ is Lipschitz continuous with Lipschitz constant $O(n^2)$
(for details, we refer to \cite[Theorem 1.1]{ESY2} or \cite[Section 5.2]{TVuniv}).

In order to show that $s_n(z)$ is close to $s_{sc}(z)$, the key observation is that $s_{sc} (z)$ can also be defined by the equation 
\begin{equation}\label{eq:s}
s_{sc}(z) = - \frac{1}{z+s_{sc}(z) }.
\end{equation}

This equation is stable, so if we can show  $s_n(z) \approx  -\frac{1}{ z+ s_n(z) }$ then it follows that 
$s_n(z) \approx s_{sc}(z)$.    This observation was due to  Bai et al.  \cite{BMT}, who used it  to prove the $n^{-1/2}$ rate of convergence  of $s_n(z)$ to $s_{sc}(z)$. 
In \cite{ESY2, ESY1, ESYweg},   Erd\H os et al. refined Bai's approach to prove local semi-circle law  at scales finer than $n^{-1/2}$, ultimately to $n^{-1} \log^C n$ \cite{ESY1}. 
Our main contribution here is to push the scale  further down to $n^{-1} \log n$, which we believe is (at most) a factor $\sqrt{\log n}$ from the truth. 

Recall that  $s_n(z)$ is the trace of $(W_n -zI)^{-1}$. By computing the diagonal entires, one can show (see  \cite[Chapter 11]{BS}, \cite{ESY1} or \cite[Lemma 39]{TVuniv}) 
\begin{equation} \label{eq:sn}
s_n(z)= \frac{1}{n} \displaystyle{\sum_{k=1}^{n} \frac{1}{-\frac{{\zeta}_{kk}}{ \sqrt{n}}-z- Y_k }},
\end{equation}
where $$Y_k=a^*_k (W_{n,k} -zI)^{-1} a_k$$ and $W_{n,k}$ is the matrix $W_n$ with the $k$-th row and column removed, and $a_k$ is the $k$-th row of $W_n$ with the $k$-th element removed. 

The entries of $a_k$ are independent of each other and of $W_{n,k}$, and have mean zero and variance $1/n$. By linearity of expectation we have 
$$\mathbf{E}(Y_k|W_{n,k})=\frac{1}{n}\tr(W_{n,k}-zI)^{-1}=(1-\frac{1}{n})s_{n,k}(z)$$ where 
$$s_{n,k}(z):= \frac{1}{n-1} {\sum_{i=1}^{n-1} \frac{1}{\lambda_i (W_{n,k}) -z}}$$ is the Stieltjes transform of $W_{n,k}$. From the Cauchy interlacing law, we can get$$\displaystyle{|{} s_n(z)- (1-\frac{1}{n}) {} s_{n,k}(z)|= O(\frac{1}{n} \int_{\mathbb{R}} \frac{1}{|x-z|^2}\,dx) =O(\frac{1}{n\eta})}=o(\delta^2)$$ and thus $$\mathbf{E}(Y_k|W_{n,k})=s_n(z)+o(\delta^2).$$

The heart of the matter now is  the following concentration result. 
\begin{lemma}
\label{YProp1} Let $M_n$ be as in Lemma \ref{lemma1}. 
For $1 \le k \le n$, $Y_k= \mathbf{E}(Y_k|W_{n,k}) +o(\delta^2)$ holds with  probability at least $1-O(n^{-C})$  for any  $z$ with $|\text{Re}(z)| \le 4$ and $\text{Im}(z) \ge {\eta}$.
\end{lemma}

To prove this lemma, we are going to make an essential use of the weighted projection lemma, as showed in the next section.  

\section{Proof  of  Lemma \ref{YProp1} and Threshold of the Local Law }\label{section:LSCLproof} 

We are going to prove Lemma \ref{YProp1} and 
the following more quantitative version of Theorem \ref{theorem:main}.

\begin{theorem} \label{theorem:LSCL}  For any constants $\epsilon, \delta, C_1>0$, there is a constant $C_2>0$ such that the following holds. 
 Let $M_n$ be a Hermitian matrix whose upper diagonal entries are independent random variables  with mean $0$ and variance $1$. 
 Assume furthermore that  for $1\le i \le n$, the vectors $X_i$, obtained by deleting the $i$-th entry of the $i$-th row vector of $M_n$, are $K$-concentrated. Then with probability at least 
$1- n^{-C_1} $, we have 
$$| N_I - n \int_I  \rho_{sc}(x) \ dx | \le \delta n \int_I  \rho_{sc}(x) \ dx , $$ 
for all interval $I \subset (-2 +\epsilon, 2-\epsilon)$ of length at least $C_2 K^2  \log n/n$. 
\end{theorem}

First, we record a lemma that provides a crude upper bound on the number of eigenvalues in short intervals.

\begin{lemma}\label{lem:crude} Let $M_n$ be a random Hermitian matrix with independent $K$-bounded entries with mean 0 and variance 1. 
For any constant $C_1>0$, there exists a constant $C_2 >0$ such that for any interval $I \subset \R$ with $|I| \ge \frac{C_2K^2 \log n}{n}$, $$N_I \ll n|I|$$ with probability at least $1-n^{-C_1}$.
\end{lemma}

This lemma is  Proposition 66 in \cite{TVuniv}, which is a variant of \cite[Theorem 5.1]{ESYweg}.
Notice that 
\begin{equation} Y_k=a^*_k (W_{n,k} -zI)^{-1} a_k=\sum_{j=1}^{n-1} \frac{|u_j (W_{n,k})^*a_k|^2}{\lambda_{j}(W_{n,k})-z}=\frac{1}{n}\sum_{j=1}^{n-1} \frac{|u_j (W_{n,k})^*X_k|^2}{\lambda_{j}(W_{n,k})-z}, \end{equation} 
where $X_k=\sqrt{n} a_k$ is the $k$-th row of $M_n$ with the $k$-th element removed. Note that the entries of $X_k$ are independent with mean 0 and variance 1.
Therefore, 
\begin{equation}\label{Ydev1}
|Y_k- \mathbf{E}(Y_k|W_{n,k})|
=\frac{1}{n}|\sum_{j=1}^{n-1} \frac{|u_j (W_{n,k})^*X_k|^2- 1}{\lambda_{j}(W_{n,k})-z} | = \frac{1}{n}|\sum_{j=1}^{n-1} \frac{R_j}{\lambda_{j}(W_{n,k})-x-\sqrt{-1}\eta } |, 
\end{equation} where $R_j:= |u_j (W_{n,k})^*X_k|^2- 1$. 
 By symmetry, we can restrict the sum to those indices $j$ where  $\lambda_j(W_{n,k})-x \ge 0$.

 Let $J$ be the set of indices $j$ such that 
  $0\le \lambda_j(W_{n,k})-x \le \eta$.   
 Since $x = \Re z, \eta =\Im z$, we have 
 \begin{eqnarray*}
&\frac{1}{n}|\sum_{j\in J} \frac{R_j}{\lambda_{j}(W_{n,k})-x-\sqrt{-1}\eta } |  \\
&\le \frac{1}{n}| \sum_{j \in J} \frac{\lambda_j(W_{n,k})-x}{(\lambda_j(W_{n,k})-x)^2+\eta^2} R_j| +\frac{1}{n}| \sum_{j \in J} \frac{\eta}{(\lambda_j(W_{n,k})-x)^2+\eta^2} R_j|\\
&\le \frac{1}{n\eta} | \sum_{j \in J} \frac{(\lambda_j(W_{n,k})-x)\eta}{(\lambda_j(W_{n,k})-x)^2+\eta^2} R_j| +\frac{1}{n\eta}| \sum_{j \in J} \frac{\eta^2}{(\lambda_j(W_{n,k})-x)^2+\eta^2} R_j|.
\end{eqnarray*}

Consider the sum $S_1:= | \sum_{j \in J} \frac{(\lambda_j(W_{n,k})-x)\eta}{(\lambda_j(W_{n,k})-x)^2+\eta^2} R_j|$. As  $0 \le \frac{(\lambda_j(W_{n,k})-x)\eta}{(\lambda_j(W_{n,k})-x)^2+\eta^2}  \le 1$, we are in position to
apply Lemma \ref{lemma:VW1}. Taking $t = C_4 K \sqrt {\log n}$  with a  sufficiently large constant $C_4$, by \eqref{ineq:VW} we have 
$$S_1 \le   \frac{C_4}{ n\eta} (2 K \sqrt{ |J| \log n } + C_4 K^2 \log n ) $$ with probability at least $1 -  C \exp (-C' C_4^2 \log n ) \ge 1 - n^{-C_4/2}$. 
By  Lemma \ref{lem:crude},  $|J| \le Bn\eta$ with probability  at least $1- n^{-C_4 }$, for some sufficiently large constant $B >0$. 
Recall $\eta:= \frac{K^2 C_3^2 \log n}{ n \delta^6} $; it follows that with probability at least $1- 2 n^{-C_4/2} $ we have 
$$S_1 \le C_4 C_3^{-1} \delta ^3(2\sqrt{B} + C_4 C_3^{-1} \delta^3). $$

Thus, for  $C_3$  sufficiently large compared to $C_4$ and $B$,  we have $S_1 \le \delta ^3$.  Similarly, we can prove the same bound for 
$S_2: = \frac{1}{n\eta}| \sum_{j \in J} \frac{\eta^2}{(\lambda_j(W_{n,k})-x)^2+\eta^2} R_j|$.

For the other eigenvalues, we divide the real line into small intervals.  For integer $l\ge 0$, let $J_l$ be the set of eigenvalues $\lambda_j(W_{n,k})$ such that $10^l \eta < \lambda_j(W_{n,k})-x \le 10^{l+1} \eta$. The number of such $J_l$ is at most $20 \log n$. By Lemma \ref{lem:crude} one has, $|J_l|\le 9B10^l n\eta$ with probability  at least $1- n^{-C_4}$, for some sufficiently large constant $B >0$. Again by Lemma \ref{lemma:VW1} (taking $t=KC_4\sqrt{\log n}$),
\begin{equation*}
\begin{split}
&\frac{1}{n}|\sum_{j\in J_l} \frac{R_j}{\lambda_{j}(W_{n,k})-x-\sqrt{-1}\eta } | \\
&\le \frac{1}{n}| \sum_{j \in J_l} \frac{\lambda_j-x}{(\lambda_j-x)^2+\eta^2} R_j| +\frac{1}{n}| \sum_{j \in J_l} \frac{\eta}{(\lambda_j-x)^2+\eta^2} R_j|\\
&\le \frac{1}{10^l n\eta}| \sum_{j \in J_l} \frac{10^l \eta(\lambda_j-x)}{(\lambda_j-x)^2+\eta^2} R_j| +\frac{1}{10^{2l} n\eta}| \sum_{j \in J_l} \frac{(10^l \eta)^2}{(\lambda_j-x)^2+\eta^2} R_j| \\
&\le \frac{2C_4 K}{10^l n \eta} (2\sqrt{|J_l|}\sqrt{\log n}+K C_4\log n)\\
&\le \delta^3 B C_4 C_3^{-1}10^{-l/2+2}
\end{split}
\end{equation*}
with probability at least $1-2C\exp(-C'C_4^2\log n)-n^{-C_4} \ge 1-n^{-C_4/2}$. 

Summing over $l$, we have $$\frac{1}{n}|\sum_{l}\sum_{j\in J_l} \frac{R_j}{\lambda_{j}(W_{n,k})-x-\sqrt{-1}\eta } | \le 200 C_4 C_3^{-1} B \delta^3 \le \delta^3$$
with probability at least $1-n^{-C_4/2 +1}$, for $C_3$ sufficiently large. This completes the proof of Lemma \ref{YProp1}.

Inserting the bounds into (\ref{eq:sn}), one has $$s_n(z)+\frac{1}{n}\sum_{k=1}^n \frac{1}{s_n(z)+z+o(\delta^2)}=0$$ with probability at least $1-O(n^{-C})$. The term $| \zeta_{kk}/\sqrt{n} | = o(\delta^2)$ as $| \zeta_{kk} | \le K$ by assumption. Comparing this equation with (\ref{eq:s}), one can use a continuity argument (see \cite{TVuniv2} for details) to obtain $|s_n(z)-s(z)|\le \delta$ with probability at least $1-O(n^{-C+100})$. 

By Lemma \ref{lemma1}, it follows that for random matrices $M_n$ with $K$-bounded entries, for any constant $C_1 >0$,  there exists a constant $C_2>0$ such that for $0 \le \delta \le1/2$ and any interval $I \subset (-3,3)$ of length at least $C_2K^2 \log n/{n\delta^8}$,
\begin{equation}\label{eq:LSCLedge}
|N_I- n \displaystyle\int_I {\rho}_{sc}(x)\,dx| \le \delta n |I|
\end{equation}
holds with probability at least $1-n^{-C_1}$. In particular,  Theorem \ref{theorem:LSCL} follows.



\section{ The infinity norm of eigenvectors } \label{section:evproof}

We prove Theorem \ref{theorem:delo} in the following more general form. 

\begin{theorem}[Optimal infinity norm of eigenvectors]  \label{theorem:delo-general} Let $M_n$ be a Hermitian matrix whose upper diagonal entries are independent  random variables with mean 0 and  variance 1. 
Further assume that for any index $1 \le i \le n$, the vector $X_i$, obtained by deleting the $i$-th entry of the $i$-th row vector of $M_n$, is $K$-concentrated. Let $W_n = \frac{1}{\sqrt{n}}M_n$. Then for 
any constant $C_1 >0$, there is a constant $C_2 >0$ such that the following holds.

\begin{itemize}
\item (Bulk case) With probability at least $1-n^{-C_1}$, for any $\epsilon >0$ and any $1\le i \le n$ with $\lambda_i(W_n) \in [-2+\epsilon, 2-\epsilon]$ there is a unit eigenvector  $u_i(W_n)$ of $\lambda_i (W_n)$ satisfying 

 $$\| u_i(W_n) \|_{\infty} \le \frac{C_2 K  \log^{1/2} n}{\sqrt n}.$$ 
 
\item (Edge case)  With probability at least $1-n^{-C_1}$, for any $\epsilon >0$ and any $1\le i \le n$ with $\lambda_i(W_n) \in [-2-\epsilon, -2+\epsilon] \cup [ 2-\epsilon, 2+\epsilon]$, there is a unit eigenvector  $u_i(W_n)$ of $\lambda_i (W_n)$ satisfying 

 $$\| u_i(W_n) \|_{\infty} \le \frac{C_2  K^2  \log n}{\sqrt n}.$$

\end{itemize}
\end{theorem}

We give here  the proof of the first part of Theorem \ref{theorem:delo-general}. The proof of the second part is somewhat different and deferred to the appendix. 
With the threshold for local semi-circle law, we are able to derive the eigenvector delocalization results thanks to the next lemma.

\begin{lemma} 
[Eq (4.3), \cite{ESY2} or Lemma 41, \cite{TVuniv}]
\label{lem:entry}
Let $$W_n=
\left(
\begin{array}{cc}
a & Y^* \\
Y & W_{n-1}
\end{array}
\right)
 $$
be an $n\times n$ Hermitian matrix for some $a\in \mathbb{C}$ and $Y \in \mathbb{C}^{n-1}$, and let 
$\left(
\begin{array}{cc}
x\\
v
\end{array}
\right)$ be an eigenvector of $W_n$ with eigenvalue $\lambda_i(W_n)$, where $x\in \mathbb{C}$ and $v \in \mathbb{C}^{n-1}$. Assume none of the eigenvalues of $W_{n-1}$ equals $\lambda_i(W_n)$. Then 
$$|x|^2= \frac{1}{1+\sum_{j=1}^{n-1} (\lambda_j (W_{n-1}) - \lambda_i(W_n))^{-2} |u_j (W_{n-1})^* Y|^2},$$
where $u_j(W_{n-1})$ is a unit eigenvector corresponding to the eigenvalue $\lambda_j (W_{n-1}).$
\end{lemma}

The assumption that the eigenvalues of $W_n$ and $W_{n-1}$ do not collide 
was taken care of in \cite[Section 3.1]{TVW}, so we can assume that the above formula makes sense in applications.

First, for the bulk case, for any $\lambda_i(W_n) \in (-2+\varepsilon, 2-\varepsilon)$, by Theorem \ref{theorem:LSCL}, one can find an interval $I \subset (-2+\varepsilon, 2-\varepsilon)$, centered at $\lambda_i(W_n)$ and with length $|I|={K^2 C \log n}/{n}$, such that $N_{I} \ge \delta_1 n|I|$ ($\delta_1 >0$ small enough)  with probability at least $1-n^{-C_1-10}$. By Cauchy interlacing law, we can find a set $J \subset \{1,\ldots,n-1\}$ with $|J| \ge N_I/2$ such that $|\lambda_j(W_{n-1})-\lambda_i(W_n)| \le |I|$ for all $j \in J$. Let $X$ be the first column of $M_n$ with the first entry removed. Then $X=\sqrt{n}Y$.

By Lemma \ref{lem:entry}, we have 
\begin{equation} \label{eq:entry}
\begin{split}
|x|^2 &=\displaystyle\frac{1}{1+\sum_{j=1}^{n-1} (\lambda_j (W_{n-1}) - \lambda_i(W_n))^{-2} |u_j (W_{n-1})^* \frac{1}{\sqrt{n}}X|^2} \\
&\le \frac{1}{1+\sum_{j \in J} (\lambda_j (W_{n-1}) - \lambda_i(W_n))^{-2} |u_j (W_{n-1})^* \frac{1}{\sqrt{n}}X|^2} \\
&\le  \frac{1}{1+ n^{-1}|I|^{-2} \sum_{j \in J}  |u_j (W_{n-1})^* X|^2} \\
&\le \frac{1}{1+100^{-1}{n^{-1}|I|^{-2}}{|J|}} \le 200 |I|/\delta_1 \le \frac{K^2 C_2^2 \log n}{ n}
\end{split}
\end{equation}
for some constant $C_2$ with  probability at least $1-n^{-C_1-10}$. The third inequality follows from \eqref{ineq:VW} by taking $t=\delta_1 K \sqrt{C \log n}$ (say).

Thus, by union bound and symmetry, $\| u_i(W_n) \|_{\infty} \le \frac{C_2 K \log^{1/2} n}{\sqrt n}$ holds with probability at least $1 -n^{-C_1}$.

\appendix
\section{Proof for the Edge case of Theorem \ref{theorem:delo-general}}\label{appendix:edge}
For the edge case in Theorem \ref{theorem:delo-general}, we use a different approach based on the next lemma.

\begin{lemma}[Interlacing identity, Lemma 37, \cite{TVuniv2}]\label{lem:interlacing} Let $W_{n-1}$ be the matrix $W_n$ with the $n$-th row and $n$-th column removed and $Y$ is the $n$-th column of $W_n$ with the $n$-th element $\zeta_{nn}/\sqrt{n}$ removed. If none of the eigenvalues of $W_{n-1}$ equals $\lambda_i(W_n)$, then
\begin{equation}
\sum_{j=1}^{n-1}\frac{|u_j(W_{n-1})^*Y|^2}{\lambda_j(W_{n-1})-\lambda_i(W_n)} =\frac{1}{\sqrt{n}} \zeta_{nn}-\lambda_i(W_n).
\end{equation}
\end{lemma}
By symmetry, it suffices to consider the case $\lambda_i(W_n) \in [2-\epsilon, 2+\epsilon]$ for $\epsilon >0$ very small. Denote $X$ the $n$-th column of $M_n$ with the $n$-th element  removed. Thus $Y=\sqrt{n}$X. By  Lemma \ref{lem:entry}, in order to show $|x|^2 \le C^4 K^4  \log^2 n/n$ (for constant $C>C_1+100$) with probability at least $1-n^{-C_1-10}$,  it is enough to show 
$$\sum_{j=1}^{n-1}\frac{|u_j (W_{n-1})^*X|^2}{(\lambda_j (M_{n-1}) - \lambda_i(M_n))^2} \ge \frac{n}{C^4 K^4 \log^2 n}.$$
By the projection lemma, $|u_j (W_{n-1})^*X| \le 10K\sqrt{C\log n} $ with  probability at least $1-10n^{-C}$.
It suffices to show that with probability at least $1-n^{-C_1-10}$, 
$$
\sum_{j=1}^{n-1}\frac{|u_j (W_{n-1})^* X|^4}{(\lambda_j (M_{n-1}) - \lambda_i(M_n))^2} \ge \frac{100n}{C^3 K^2 \log n}.
$$

By Cauchy-Schwardz inequality, it is enough to show for some integers $1\le T_{-} < T_{+} \le n-1$ that
$$
\sum_{T_- \le j \le T_+}\frac{|u_j (W_{n-1})^* Y|^2}{|\lambda_j (W_{n-1}) - \lambda_i(W_n)|} \ge \frac{10\sqrt{T_+ - T_-}} {C^{1.5} K \sqrt{\log n}}.
$$

By Lemma \ref{lem:interlacing}, we are going to show for some integers $T_{+}, T_{-}$ satisfying $T_{+} - T_{-} = O(\log n)$ (the choice of $T_{+}, T_{-}$ will be given later) that 
\begin{equation}\label{eq:other}
|\sum_{j \ge T_+ \text{or} j \le T_-} \frac{|u_j (W_{n-1})^* Y|^2}{\lambda_j (W_{n-1}) - \lambda_i(W_n)}| \le 2-\epsilon-\frac{10\sqrt{T_+ - T_-}} {C^{1.5} K \sqrt{\log n}}+o(1),
\end{equation}
with probability at least $1-n^{-C_1-10}$.

Let $\eta =\frac{K^2 C \log n}{n\delta^8}$  with constant $\delta= \epsilon/1000$. Divide the real line into disjoint intervals $I_k$ for $k \ge 0$ where $I_0 = (\lambda_i(W_n)- \eta, \lambda_i(W_n)+ \eta)$. For $1 \le k \le k_0=\log^{0.9} n$ (say),  $|I_k|$ has length $2\eta \delta^{-8k} = o(1)$ and  
$$
I_k = (\lambda_i(W_n)-\beta_k \eta, \lambda_i(W_n)-\beta_{k-1} \eta ] \cup [\lambda_i(W_n) + \beta_{k-1} \eta, \lambda_i(W_n)+ \beta_k \eta ),
$$
where we denote by $\beta_k = \sum_{s=0}^{k} \delta^{-8s}$. The distance from $\lambda_i(W_n)$ to the interval $I_k$ satisfies 
$$
\mathrm{dist}(\lambda_i(W_n), I_k) \ge \beta_{k-1} \eta.
$$

For each such interval, by \eqref{eq:LSCLedge}, for sufficiently large constant $C>0$, the number of eigenvalues $|J_k|=N_{I_k}\le n\alpha_{I_k} |I_k| +\delta^{k+1} n |I_k| $ with probability at least $1-n^{-C_1-100}$, where $\alpha_{I_k}=\int_{I_k} \rho_{sc}(x) dx/ |I_k|$.

For the $k$-th interval, by  \eqref{ineq:VW} taking $t=K\sqrt{C \log n}$,  we have that, with probability at least $1-C''\exp(-C'C\log n) \ge 1-n^{-C_1-100}$ for sufficiently large $C$,
\begin{equation*}
\begin{split}
&\frac{1}{n} \sum_{j\in J_k} \frac{|u_j(W_{n-1})^*X|^2}{|\lambda_j (W_{n-1})-\lambda_i(W_n)|} \le \frac{1}{n}\frac{1}{\text{dist}(\lambda_i(W_n), I_k)}\sum_{j\in J_k}|u_j(W_{n-1})^*X|^2 \\
&\le \frac{1}{n}\frac{1}{\text{dist}(\lambda_i(W_n), I_k)}(|J_k|+K\sqrt{|J_k|}\sqrt{C\log n}+C K^2  \log n)\\
&\le \frac{\alpha_{I_k} |I_k|}{\text{dist}(\lambda_i(W_n), I_k)}+\frac{1}{n\eta \beta_{k-1}} (n\delta^{k+1}|I_k| + \sqrt{2}K \sqrt{C\log n} \sqrt{n|I_k|} + CK^2 \log n) \\
&\le \frac{\alpha_{I_k} |I_k|}{\text{dist}(\lambda_i(W_n), I_k)}+10\delta^{k-7}.
\end{split}
\end{equation*}

For $k \ge k_0 + 1$, let the interval $I_k$'s have the same length of $| I _{k_0}|= 2 \delta^{-8k_0} \eta$. Note that the number of such intervals is bounded crudely by $o(n)$.  By \eqref{eq:LSCLedge}, the number of eigenvalues $|J_k| \le n\alpha_{I_k} |I_k| +\delta^{k_0+1} n |I_k| $ with probability at least $1-n^{-C_1-100}$. And the distance from $\lambda_i(W_n)$ to the interval $I_k$ satisfies 
$$
\mathrm{dist}(\lambda_i(W_n), I_k) \ge \beta_{k_0-1} \eta + (k-k_0) |I_{k_0}|.
$$

The contribution of such intervals can be computed similarly 
\begin{equation*}
\begin{split}
&\frac{1}{n} \sum_{j\in J_k} \frac{|u_j(W_{n-1})^*X|^2}{|\lambda_j (W_{n-1})-\lambda_i(W_n)|} \le \frac{1}{n}\frac{1}{\text{dist}(\lambda_i(W_n), I_k)}\sum_{j\in J_k}|u_j(W_{n-1})^*X|^2 \\
&\le \frac{1}{n}\frac{1}{\text{dist}(\lambda_i(W_n), I_k)}(|J_k|+K\sqrt{|J_k|}\sqrt{C\log n}+C K^2  \log n)\\
&\le \frac{\alpha_{I_k} |I_k|}{\text{dist}(\lambda_i(W_n), I_k)}+\frac{ \delta^{k_0} }{ k-k_0}
\end{split}
\end{equation*}
with probability at least $1-n^{-C_1-100}$.
 
Summing over all intervals for $k\ge 10$ (say), we obtain
\begin{equation}\label{eq:upperbd1}
|\sum_{j \ge T_+ \text{or} j \le T_-} \frac{|u_j (W_{n-1})^* Y|^2}{\lambda_j (W_{n-1}) - \lambda_i(W_n)}| \le |\sum_{I_k}\frac{\alpha_{I_k} |I_k|}{\text{dist}(\lambda_i(W_n), I_k)}|+\delta.
\end{equation}

On the other hand, it follows from Riemann integration of the principal value integral that
$$\sum_{I_k}\frac{\alpha_{I_k} |I_k|}{\text{dist}(\lambda_i(W_n), I_k)}=p.v.\int_{-2}^2 \frac{\rho_{sc}(x)}{\lambda_i(W_n) - x} \ dx +o(1),$$
where $$p.v.\int_{-2}^2 \frac{\rho_{sc}(x)}{\lambda_i(W_n) - x} \ dx :=\lim_{\varepsilon\rightarrow 0} \int_{-2\le x\le 2, |x-\lambda_i(W_n)| \ge \varepsilon} \frac{\rho_{sc} (x)}{ \lambda_i(W_n)-x} \ dx.$$

From the explicit formula for the Stieltjes transform and from residue calculus, one obtains
$$p.v.\int_{-2}^2 \frac{\rho_{sc}(x)}{x-\lambda_i(W_n)} \ dx =-\lambda_i(W_n)/2$$ for $|\lambda_i(W_n)| \le 2$, and with the right-hand side replaced by $-\lambda_i(W_n)/2 + \sqrt{\lambda_i(W_n)^2-4}/2$ for $|\lambda_i(W_n)| > 2$. Finally, we always have 
\begin{equation}\label{eq:mainpart}
|\sum_{I_k}\frac{\alpha_{I_k} |I_k|}{\text{dist}(\lambda_i(W_n), I_k)}| \le 1+2\epsilon.
\end{equation}

Now for the rest of eigenvalues  that satisfy $|\lambda_i(W_n)-\lambda_j(W_{n-1})|\le |I_0|+|I_1|+\ldots+|I_{10}|\le 4\eta/\delta^{80}$, by Theorem \ref{theorem:LSCL} and Cauchy interlacing law, the number of eigenvalues is at most ${T_+-T_-} \le 8n\eta/\delta^{80} =8C K^2 \log n/\delta^{88}$ with probability at least $1-n^{-C_1-100}$ for sufficiently large constant $C>0$. Thus 
\begin{equation}\label{eq:upperbd2}
\frac{\sqrt{T_+ - T_-}}{C^{1.5} K \sqrt{\log n}} \le \frac{10}{ \delta^{44} C} \le \epsilon/1000,
\end{equation}
by choosing $C$ sufficiently large compared to $\delta^{-44}$. Thus, from \eqref{eq:other}, \eqref{eq:upperbd1}, \eqref{eq:mainpart} and \eqref{eq:upperbd2}, we have proved that there exits a constant $C>0$ such that with probability at least $1-n^{-C_1-10}$,
$$|x| \le \frac{C^2 K^2 \log n}{\sqrt{n}}.$$

The conclusion of the second part of Theorem \ref{theorem:delo} follows from symmetry and union bounds.

\section{Local Marchenko-Pastur law for random covariance matrix and delocalization of singular vectors}\label{appendix:covariance}

In this appendix, we extend the results obtained for random Hermitian matrices discussed in the previous sections to random covariance matrices, focusing on the changes needed for the proofs. Interested reader can refer to closely related papers  \cite{TVcov} and \cite{KWcov} (see also \cite{ESYY, PY}).

Let $M=M_{p,n} =(\zeta_{ij})_{1\le i\le p, 1\le j \le n}$ be a $p\times n$ matrix, where $p=p(n)$ is an integer such that $p\le n$ and $\lim_{n \rightarrow \infty} p/n=y\in (0,1]$. Assume the entries of $M_{n,p}$ are independent random variables with mean zero and variance one. For such a $p\times n$ random matrix $M$, we form the $n \times n$ (sample) covariance matrix $ W= W_{p,n}= \frac{1}{n} M^*M$. This (non-negative definite) matrix has at most $p$ non-zero eigenvalues which are ordered as 
$$
0 \le {\lambda}_1 {(W)} \le {\lambda}_2 {(W)} \le \ldots \le {\lambda}_p {(W)}.
$$ 

Denote by $\sigma_1(M),\ldots, \sigma_p(M)$ the singular values of $M$. It is easy to see that $\sigma_i(M) = \sqrt{n} \lambda_i (W)^{1/2}$. From the singular value decomposition, there exist orthonormal bases $\{u_1,\ldots,u_p \}$ for $\mathbb{C}^n$ and $\{v_1,\ldots,v_p\}$ for $\mathbb{C}^p$ such that 
$Mu_i=\sigma_i v_i$ and $M^*v_i=\sigma_i u_i.$

A fundamental result concerning the asymptotic limiting behavior of ESD for large covariance matrices is the Marchenko--Pastur Law (see \cite{BS1} and  \cite{MP}). 
\begin{theorem}\label{thm:MP}
(Marchenko--Pastur  Law)
Assume the entries of matrix $M \in \C^{p\times n}$ are independent random variables with mean zero and variance one and  $\lim_{n \rightarrow \infty} p/n=y \in (0,1]$. Then the empirical spectral distribution of the matrix $W= \frac{1}{n} M^*M$ converges with probability 1 to the Marchenko-Pastur  Law with a density function  $${{\rho}}_{MP,y}(x) := \frac{1}{2 \pi xy} \sqrt{(b-x)(x-a)} \mathbf{1}_{[a,b]}(x),$$ where $$a:=(1-\sqrt{y})^2, b:=(1+\sqrt{y})^2.$$
\end{theorem}

The \emph{hard edge} of the limiting support of spectrum refers to the left edge $a$ when $y=1$ where it gives rise to a singularity of $x^{-1/2}$. The cases of left edge $a$ when $y <1$ and the right edge $b$ regardless of the value of $y$ are usually called the \emph{soft edge}. Recent progress on studying the local convergence to Marchenko--Pastur Law includes \cite{ESYY, PY, TVcov, KWcov} for the soft edge and \cite{CMS, TVsmall} for the hard edge. We focus on improving the previous results for the soft edge in this appendix. 

Our main results for the random covariance matrices are the following local Marchenko--Pastur law (LMPL) and the delocalization property of singular vectors.

\begin{theorem} \label{thm:LMPL}
For any constants $\epsilon, \delta, C_1 >0$, there exists a constant $C_2>0$ such that the following holds. Assume $\lim_{n\rightarrow \infty}p/n= y$ for some $0<y\le 1$. Let $M=M_{p,n}=(\zeta_{ij})_{1\le i\le p, 1\le j \le n}$ be a random matrix whose entries are independent $K$-bounded random variables with mean 0 and variance 1. Consider the covariance matrix $W= \frac{1}{n} M^*M$. Then with probability at least $1-n^{-C_1}$, one has
\[
|N_I (W_{n,p}) -p \int_I {{}\rho}_{MP,y}(x)\,dx| \le {\delta} p  \int_I {{}\rho}_{MP,y}(x)\,dx.
\]
for any interval $I \subset (a+\epsilon, b-\epsilon)$ of length at least $C_2 K^2 \log n/n.$
\end{theorem}

\begin{theorem}[Delocalization of singular vectors]  \label{thm:singularvector} Let $M_{p,n}$ be as in Theorem \ref{thm:LMPL}. For any constant $C_1 >0$, there is a constant $C_2 >0$ such that the following holds.
\begin{itemize}
\item (Bulk case) With probability 
at least $1- n^{-C_1}$, for any $\epsilon >0$ and any $1\le i \le p$ such that $\sigma_i(M_{p,n})^2/n \in [a+\epsilon, b-\epsilon]$, there is a left singular vector   $u_i$ corresponding to 
$\sigma_i (M_{p,n} )$ such that  $$\| u_i \|_{\infty} \le \frac{C_2 K \log^{1/2} n}{\sqrt n}.$$  
The same holds for right singular vectors.  
 
\item (Edge case) With probability at least $1-n^{-C_1}$, for any
 $\epsilon >0$ and any $1\le i \le p$ such that $\sigma_i(M_{p,n})^2/n \in [a-\epsilon, a+\epsilon] \cup [ b-\epsilon, b+\epsilon]$ if $a\neq 0$ and $\sigma_i(M_{p,n})^2/n \in  [ 4-\epsilon, 4]$ if $a=0$,
 there is a left singular vector $u_i$ corresponding to 
$\sigma_i (M_{n,p} )$ such that
  $$\| u_i \|_{\infty} \le \frac{C_2 K^2 \log n}{\sqrt n}.$$ The same holds for right singular vectors. 

\end{itemize}
\end{theorem} 

\remark Theorem \ref{thm:LMPL} and Theorem \ref{thm:singularvector} actually hold for a larger class of matrices, using the $K$-concentration introduced in the previous sections. For instance, Theorem \ref{thm:LMPL} holds for random matrices $M_{p,n}=(\zeta_{ij})$ whose entries are independent random variables with mean 0 and variance 1, and the row vectors are $K$-concentrated. And Theorem \ref{thm:singularvector} holds if we further assume the column vectors of $M_{p,n}$ are also $K$-concentrated. Indeed, the $K$-bounded assumption is only used to guarantee $K$-concentration.

\subsection{Proof of Theorem \ref{thm:LMPL}}

Similarly to the Hermitian case, we compare the Stieltjes transform of  $W$
$$s(z):=\frac{1}{p} \sum_{i=1}^{p} \frac{1}{\lambda_i(W)-z}$$ 
with that of the Marchenko--Pastur Law
$$
s_{MP,y}(z):= \int_{\R} \frac{1}{x-z} \rho_{MP,y}(x)\,dx=\int_{a}^{b} \frac{1}{2\pi xy(x-z)} \sqrt{(b-x)(x-a)}\,dx.
$$ 

The explicit expression of $s_{MP,y}(z)$ is given by (see \cite{BS})
$$s_{MP,y}(z)=-\frac{y+z-1-\sqrt{(y+z-1)^2-4yz}}{2yz},$$ where we take the branch of $\sqrt{(y+z-1)^2-4yz}$ with cut at $[a,b]$ that is asymptotically $y+z-1$ as $z$ tends to infinity. Note that  it is uniquely defined by the equation 
$$
s_{MP,y}(z)+\frac{1}{y+z-1+yzs_{MP,y}(z)}=0.
$$

We will show that $s(z)$ satisfies a similar equation.

The analogue of Lemma \ref{lemma1} is the following lemma.

\begin{lemma}
\emph{(Lemma 29, \cite{TVcov})}  
\label{ESDStie}
Let $M_{n,p}$ be a random matrix with independent $K$-bounded entries with mean 0 and variance 1. Assume $\lim_{n\rightarrow +\infty}p/n=y \in (0,1]$. Let $1/n< \eta < 1/10$, and $L_1, L_2, \varepsilon, \delta >0$. For any constant $C_1 >0$, there exists a constant $C>0$ such that if one has the bound $$|s(z)-s_{MP,y}(z)| \le \delta$$ with (uniformly) probability at least $1-n^{-C}$ for all $z$ with $ L_1 \le \text{Re}(z) \le L_2$ and $\text{Im}(z) \ge \eta$. Then for any interval $I$ in $[L_1-\varepsilon, L_2+\varepsilon]$ with $|I| \ge \text{max}(2\eta, \frac{\eta}{\delta} \log \frac{1}{\delta})$, one has $$|N_I- p \int_I {\rho}_{MP,y}(x)\,dx| \le \delta p |I|$$ with probability at least $1-n^{-C_1}$. 
\end{lemma}

The objective is to show 
\begin{equation} \label{eq:diff}
|s(z)-s_{MP,y}(z)| =o (\delta)
\end{equation}
with probability at least $1-n^{-C}$ for any $z$ in the region $R_y$, where $$R_y = \{ z \in \mathbb{C}: |z| \le 10, a-\epsilon \le \text{Re}(z) \le b+\epsilon, \text{Im}(z) \ge {\eta} \}$$ if $y \neq 1$, and $$R_y = \{ z \in \mathbb{C}: |z| \le 10, \epsilon \le \text{Re}(z) \le 4+\epsilon, \text{Im}(z) \ge {\eta} \}$$ if $y=1$. We use the parameter 
$$
\eta:=\frac{K^2 C^2 \log n }{n \delta^6},
$$ 
where $C=C_1+10^4$. Note that in the defined region $R_y$, $| s_{MP,y}(z) | = O(1)$.


First, by Schur's complement, one can rewrite
\begin{equation} \label{eq:1.6}
s(z)= \frac{1}{p} \tr(W^*-z I)^{-1}=\frac{1}{p} \displaystyle{\sum_{k=1}^{p} \frac{1}{\xi_{kk}-z- Y_k }}
\end{equation}
where $Y_k=a^*_k (W_{k} -zI)^{-1} a_k$, and $W_{k}$ is the matrix $W^*=\frac{1}{n}M M^*=(\xi_{ij})_{1\le i,j \le p}$ with the $k$-th row and $k$-th column removed, and $a_k$ is the $k$-th row of $W$ with the $k$-th element removed. Let $M_k$ be the $(p-1)\times n$ minor of $M$ with the $k$-th row removed and $X_i^* \in \mathbb{C}^n (1\le i \le p)$ be the rows of $M$. Thus $\xi_{kk}={X_k}^*X_k/n= \|X_k\|^2/n, a_k=\frac{1}{n}M_k X_{k}, W_k=\frac{1}{n}M_k M_k^*$. Thus
\begin{equation*}
Y_k  = \sum_{j=1}^{p-1} \frac{|a_k ^* v_j(M_k)|^2}{\lambda_j(W_k)-z}
= \sum_{j=1}^{p-1} \frac{1}{n} \frac{\lambda_j(W_k)|X_k^* u_j(M_k)|^2}{\lambda_j(W_k) -z}
\end{equation*}
where $u_1(M_k),\ldots,u_{p-1}(M_k) \in \mathbb{C}^{n}$ and $v_1(M_k),\ldots,v_n(M_k) \in \mathbb{C}^{p-1}$ are orthonormal right and left singular vectors of $M_k$. Here we use the fact that $a_k^* v_j(M_k)=\frac{1}{n} X_k^* M_k^*v_j(M_k)=\frac{1}{n}\sigma_j(M_k)X_k^*u_j(M_k)$ and $\sigma_j(M_k)^2=n\lambda_j(W_k)$.

The entries of $X_k$ are independent of each other and of $W_{k}$, and have mean $0$ and variance $1$. Since $u_j(M_k)$ are unit vectors, by linearity of expectation we have 
$$
\E(Y_k|W_{k})=\sum_{j=1}^{p-1}\frac{1}{n}\frac{\lambda_j(W_k)}{\lambda_j(W_k) -z}
=\frac{p-1}{n}+\frac{z}{n} \sum_{j=1}^{p-1}\frac{1}{\lambda_j(W_k)-z}=\frac{p-1}{n}(1+zs_k(z)),
$$ 
where 
$$
s_{k}(z)= \frac{1}{p-1} \displaystyle{\sum_{i=1}^{p-1} \frac{1}{\lambda_i (W_{k}) -z}}
$$ 
is the Stieltjes transform of  $W_{k}$. By Cauchy interlacing law, we have
$$
| s(z)- (1-\frac{1}{p})  s_{k}(z)|= O(\frac{1}{p} \int_{\mathbb{R}} \frac{1}{|x-z|^2}\,dx) =O(\frac{1}{p\eta}).
$$ 

Thus 
$$
\E(Y_k|W_{k})=\frac{p-1}{n}+z\frac{p}{n}s(z)+O(\frac{1}{n\eta})=\frac{p-1}{n}+z\frac{p}{n}s(z)+o(\delta^2).$$

On the other hand,  $Y_k$ is concentrated about $\E(Y_k|W_{k})$ with high probability:
\begin{lemma}
\label{YProp}
Let $M_{n,p}$ be as in Lemma \ref{ESDStie}. For $1 \le k \le p$, $Y_k= \mathbf{E}(Y_k|W_{k}) +o(\delta^2)$ holds with probability at least $1-O(n^{-C})$ for any $z$ in the region $R_y$.
\end{lemma}

 To prove Lemma \ref{YProp}, we estimate   
 \begin{equation}\label{Ydev2}
 \begin{split}
 Y_k-\mathbf{E}(Y_k|W_{k})
 &=\displaystyle\sum_{j=1}^{p-1} \frac{\lambda_j(W_k)}{n} \left( \frac{|X_k^* u_j(M_k)|^2 -1}{\lambda_j(W_k)-z} \right) =\frac{1}{n} \sum_{j=1}^{p-1} \frac{\lambda_j(W_k)}{\lambda_j(W_k)-x-\sqrt{-1} \eta }R_j
 \end{split}
 \end{equation}
where $R_j =|X_k^* u_j(M_k)|^2 -1 $. Note that $\lambda_j(W_k) = O(1)$. The estimation of (\ref{Ydev2}) is a repetition of the calculation in (\ref{Ydev1}). Interested reader are encouraged to work out the details.
Inserting the bounds to (\ref{eq:1.6}), we have $$s(z)+ \frac{1}{y+z-1+yz s(z) +o(\delta^2)}=0$$with probability at least $1-O(n^{-C})$. By a continuity argument (see for instance \cite{KWcov}), one has $|s(z)- s_{MP,y}(z)|=o(\delta)$ with probability at least $1-n^{-C+100}$ (say).  By Lemma \ref{ESDStie}, we have showed that for any constants $\epsilon, C_1 >0$, there exists a constant $C_2 >0$ such that for $0<\delta<1/2$ and any interval $I \subset (a-\epsilon,b+\epsilon)$ if $a\neq 0$ or $I \subset (\epsilon,4+\epsilon)$ if $a=0$ of length at least $C_2 K^2 \log n/n\delta^8,$ with probability at least $1-n^{-C_1}$, 
\begin{equation}\label{eq:LMPLedge}
|N_I -p \displaystyle\int_I {{}\rho}_{MP,y}(x)\,dx| \le {\delta} p |I|.
\end{equation}

In particular, Theorem \ref{thm:LMPL} follows.


\subsection{Proof of Theorem \ref{thm:singularvector}}
To prove the delocalization of singular vectors, we need the following formula to express the entries of a singular vector in terms of the singular values and singular vectors of a minor. It is enough to prove the delocalization for the right (unit) singular vectors.
\begin{lemma} [Corollary 25, \cite{TVcov}]\label{lem:singularco}
 Let $p,n \ge 1$, and let 
$$M_{p,n}=
\left(
\begin{array}{cc}
M_{p,n-1} & X
\end{array}
\right)
 $$
be a $p\times n$ matrix for some $X\in \mathbb{C}^p$, and let 
$\left(
\begin{array}{cc}
u\\
x
\end{array}
\right)$ be a right unit singular vector of $M_{p,n}$ with singular value $\sigma_i(M_{p,n})$, where $x\in \mathbb{C}$ and $u \in \mathbb{C}^{n-1}$. Suppose that none of the singular values of $M_{p,n-1}$ are equal to $\sigma_i(M_{p,n})$. Then 
$$|x|^2= \frac{1}{1+\sum_{j=1}^{\text{min}{(p,n-1)}} \frac{\sigma_j(M_{p,n-1})^2}{(\sigma_j(M_{p,n-1})^2 - \sigma_i(M_{p,n})^2)^2} |v_j (M_{p,n-1})^* X|^2},$$
where $v_1(M_{p,n-1}), \ldots, v_{\text{min}{(p,n-1)}}(M_{p,n-1}) \in \mathbb{C}^p$ is an orthonormal system of left singular vectors corresponding to the non-trivial singular values of $M_{p,n-1}.$

In a similar vein, if $$M_{p,n}=
\left(
\begin{array}{cc}
M_{p-1,n}\\
Y^*
\end{array}
\right)$$
for some $Y\in \mathbb{C}^n$, and $\left(
\begin{array}{cc}
v\\
y
\end{array}
\right)$ is a left unit singular vector of $M_{p,n}$ with singular value $\sigma_i(M_{p,n})$, where $y\in \mathbb{C}$ and $v\in \mathbb{C}^{p-1}$, and none of the singular values of $M_{p-1,n}$ are equal to $\sigma_i(M_{p,n})$. Then 
$$|y|^2= \frac{1}{1+\sum_{j=1}^{\text{min}{(p-1,n)}} \frac{\sigma_j(M_{p-1,n})^2}{(\sigma_j(M_{p-1,n})^2 - \sigma_i(M_{p,n})^2)^2} |u_j (M_{p-1,n})^* Y|^2},$$
where $u_1(M_{p-1,n}), \ldots, u_{\text{min}{(p-1,n)}}(M_{p-1,n}) \in \mathbb{C}^n$ is an orthonormal system of right singular vectors corresponding to the non-trivial singular values of $M_{p-1,n}.$
\end{lemma}

First, if $\lambda_i (W_{p,n})$ lies within the bulk of spectrum, by Theorem \ref{thm:LMPL}, one can find an interval $I \subset (a+\varepsilon, b-\varepsilon)$, centered at $\lambda_i (W_{p,n})$ and with length $|I|={K^2 C \log n}/n$ such that $N_{I} \ge \delta_1 n|I|$ ($\delta_1 >0$ small constant) with probability at least $1-n^{-C_1-10}$. By Cauchy interlacing law, we can find a set $J \subset \{1,\ldots,p\}$ with $|J| \ge N_I/2$ such that $|\lambda_j(W_{p,n-1})-\lambda_i(W_{p,n})| \le |I|$ for all $j \in J$. Thus
\begin{equation*} 
\begin{split}
&\sum_{j=1}^{\text{min}{(p,n-1)}} \frac{\sigma_j(M_{p,n-1})^2 }{(\sigma_j (M_{p,n-1})^2 - \sigma_i(M_{p,n} )^2 )^2 } |v_j (M_{p,n-1})^* X|^2 \\
&\ge \frac{1}{n}  \sum_{j \in J} \frac{ \lambda_j(W_{p,n-1})}{(\lambda_j (W_{p,n-1}) - \lambda_i(W_{p,n}))^2} |v_j (M_{p,n-1})^* X|^2\\
&\ge n^{-1}|I|^{-2} \sum_{j \in J}  |v_j (M_{p,n-1})^* X|^2 \ge 100^{-1}n^{-1}|I|^{-2} |J| \\
&\ge \frac{\delta_1}{200} |I|^{-1} \ge \frac{n}{K^2 C_2^2 \log n}
\end{split}
\end{equation*}
with  probability at least $1-n^{-C_1-10}$ for some constant $C_2>0$. The fourth inequality follows from \eqref{ineq:VW} by taking $t=\delta_1 K \sqrt{C} \log n$.

Thus, by Lemma \ref{lem:singularco} and the union bound, $|x| \le \frac{C_2 K \log^{1/2} n}{\sqrt n}$ with probability at least $1 -n^{-C_1-10}$. By symmetry and union bounds, $\| u_i(M_{p,n}) \|_{\infty} \le \frac{C_2 K \log^{1/2} n}{\sqrt n}$ holds with probability at least $1 -n^{-C_1}$.

For the edge case, we consider $| \lambda_i(W_{p,n})-a | =o(1)$ $(a\neq 0)$ or $|\lambda_i(W_{p,n}) - b| =o(1)$. We first record  an analogue of Lemma \ref{lem:interlacing}.

\begin{lemma} [Interlacing identity for singular values, Lemma 3.5 \cite{KWcov}] \label{lem:cauchyid}
Assume the notations in Lemma \ref{lem:singularco}, then for every $i$,
\begin{equation}\label{eq:lace1}
\sum_{j=1}^{{\text{min}{(p,n-1)}}} \frac{\sigma_j(M_{p,n-1})^2 |v_j (M_{p,n-1})^* X|^2}{\sigma_j(M_{p,n-1})^2 - \sigma_i(M_{p,n})^2} = ||X||^2 - \sigma_i(M_{p,n})^2.
\end{equation}
Similarly, we have 
\begin{equation}\label{eq:lace2}
\sum_{j=1}^{{\text{min}{(p-1,n)}}} \frac{\sigma_j(M_{p-1,n})^2 |u_j (M_{p-1,n})^* Y|^2}{\sigma_j(M_{p-1,n})^2 - \sigma_i(M_{p,n})^2} = ||Y||^2 - \sigma_i(M_{p,n})^2.
\end{equation}
\end{lemma}





By the union bound and Lemma \ref{lem:singularco}, in order to show $|x|^2 \le C^4 K^4  \log^2 n/n$ with probability at least $1-n^{-C_1-10}$ for some large constant $C>C_1+100$,  it is enough to show 
$$
\sum_{j=1}^{\mathrm{min}(p,n-1)}\frac{\sigma_j(M_{p,n-1})^2}{(\sigma_j(M_{p,n-1})^2 - \sigma_i(M_{p,n})^2 )^2} |v_j (M_{p,n-1})^*X|^2 \ge \frac{n}{C^4 K^4  \log^2 n}.
$$

By the projection lemma, $|v_j (M_{p,n-1})^*X| \le 10K\sqrt{C\log n} $ with  probability at least $1-10n^{-C}$. It suffices to show that with probability at least $1-n^{-C_1-10}$, 
$$
\sum_{j=1}^{\mathrm{min}(p,n-1)}\frac{\sigma_j(M_{p,n-1})^2}{(\sigma_j(M_{p,n-1})^2 - \sigma_i(M_{p,n})^2 )^2} |v_j (M_{p,n-1})^*X|^4 \ge \frac{100n}{C^3 K^2 \log n}.
$$

By Cauchy-Schwardz inequality and note that $|\sigma_i(M_{p,n-1})| \le 10\sqrt{n}$ almost surely (See \cite{Geman, YBK}), it is enough to show for some integers $1\le T_{-} < T_{+} \le \text{min} (p,n-1)$ (the choice of $T_{-}, T_{+}$ will be given later),
$$
\sum_{T_- \le j \le T_+}\frac{ \frac{1}{n} \sigma_j(M_{p,n-1})^2}{|\sigma_j(M_{p,n-1})^2 - \sigma_i(M_{p,n})^2 |} |v_j (M_{p,n-1})^*X|^2 \ge \frac{100\sqrt{T_+ - T_-}} {C^{1.5} K \sqrt{\log n}}.
$$

On the other hand, by the projection lemma, with probability at least $1- n^{-C_1-100}$, $\| X \|^2/n = y + o(1)$. By (\ref{eq:lace1}) in Lemma \ref{lem:cauchyid}, 
\begin{equation}\label{total}
 \sum_{j=1}^{{\text{min}{(p,n-1)}}} \frac{1}{n}\frac{\sigma_j(M_{p,n-1})^2 |v_j (M_{p,n-1})^* X|^2}{\sigma_j(M_{p,n-1})^2 - \sigma_i(M_{p,n})^2} =y+o(1)-\lambda_i(W_{p,n}).
\end{equation}

It is enough to evaluate
\begin{equation}\label{partial}
\begin{split}
&\sum_{j \ge T_+ \mathrm{ or } j \le T_-} \frac{\lambda_j(W_{p,n-1}) |v_j (M_{p,n-1})^* X|^2}{ \lambda_j(W_{p,n-1}) - \lambda_i(W_{p,n})}.  
\end{split}
\end{equation}

The estimation of \eqref{partial} is similar to that of \eqref{eq:other}. We divide the real line into disjoint intervals $I_k$ for $k \ge 0$. Let $\eta=\frac{K^2 C \log n}{n\delta^8}$  with small constant $\delta\le 0.01$.  Denote $\beta_k = \sum_{s=0}^{k} \delta^{-8s}$.  Let $I_0 = (\lambda_i(W_{p,n})- \eta, \lambda_i(W_{p,n})+ \eta)$. For $1 \le k \le k_0=\log^{0.9} n$ (say), 
$$
I_k = (\lambda_i(W_{p,n})-\beta_k \eta, \lambda_i(W_{p,n})-\beta_{k-1} \eta] \cup [\lambda_i(W_{p,n}) + \beta_{k-1} \eta, \lambda_i(W_{p,n})+ \beta_k \eta ),
$$
thus $| I_k |= 2 \delta^{-8k} \eta = o(1)$ and the distance from $\lambda_i(W_{p,n})$ to the interval $I_k$ satisfies 
$$
\mathrm{dist}(\lambda_i(W_{p,n}), I_k) \ge \beta_{k-1}\eta.
$$

For each such interval, by \eqref{eq:LMPLedge}, for sufficiently large constant $C>0$, the number of eigenvalues $|J_k|=N_{I_k}\le p\alpha_{I_k} |I_k| +\delta^{k+1} p |I_k| $ with probability at least $1-n^{-C_1-100}$, where $\alpha_{I_k}=\int_{I_k} \rho_{MP,y}(x) dx/ |I_k|$.

Taking $t=K\sqrt{C\log n}$ in \eqref{ineq:VW} for $C$ sufficiently large, it follows that with probability at least $1-C''\exp(-C'C\log n) \ge 1-n^{-C_1-100}$,
\begin{equation*}
\begin{split}
&\frac{1}{n} \sum_{j\in J_k} \frac{|\lambda_j(W_{p,n-1})| |v_j(M_{p,n-1})^*X|^2}{|\lambda_j (W_{p,n-1})-\lambda_i(W_{p,n})|} \le \frac{1}{n}(1+\frac{\lambda_i(W_{p,n})}{\text{dist}(\lambda_i(W_n), I_k)} ) \sum_{j\in J_k}|v_j(W_{p,n-1})^*X|^2 \\
&\le \frac{1}{n}(1+\frac{\lambda_i(W_{p,n})}{\text{dist}(\lambda_i(W_{p,n}), I_k)} ) (|J_k|+K\sqrt{|J_k|}\sqrt{C\log n}+C K^2  \log n)\\
&\le \frac{1}{n}(1+\frac{\lambda_i(W_{p,n})}{\text{dist}(\lambda_i(W_{p,n}), I_k)} ) (p\alpha_{I_k} |I_k| + \delta^k p |I_k| + \sqrt{2}K \sqrt{C\log n} \sqrt{n} \sqrt{| I_k |}+C K^2 \log n)\\
&\le y (1+\frac{\lambda_i(W_{p,n})}{\text{dist}(\lambda_i(W_{p,n}), I_k)} ) \alpha_{I_k} |I_k| + 100 \delta^{k-7}.
\end{split}
\end{equation*}

For $  k \ge k_0 + 1$, let the intervals $I_k$'s have the same length of $| I _{k_0}|= 2 \delta^{-8k_0} \eta$. Note that the number of such intervals is bounded crudely by $o(n)$. The distance from $\lambda_i(W_{p,n})$ to the interval $I_k$ satisfies 
$$
\mathrm{dist}(\lambda_i(W_{p,n}), I_k) \ge \beta_{k_0-1} \eta + (k-k_0) |I_{k_0}|.
$$

The contribution of such intervals can be estimated similarly by
\begin{equation*}
\begin{split}
&\frac{1}{n} \sum_{j\in J_k} \frac{|\lambda_j(W_{p,n-1})| |v_j(M_{p,n-1})^*X|^2}{|\lambda_j (W_{p,n-1})-\lambda_i(W_{p,n})|} \le \frac{1}{n}(1+\frac{\lambda_i(W_{p,n})}{\text{dist}(\lambda_i(W_n), I_k)} ) \sum_{j\in J_k}|v_j(W_{p,n-1})^*X|^2 \\
&\le \frac{1}{n}(1+\frac{\lambda_i(W_{p,n})}{\text{dist}(\lambda_i(W_{p,n}), I_k)} ) (|J_k|+K\sqrt{|J_k|}\sqrt{C\log n}+C K^2  \log n)\\
&\le y (1+\frac{\lambda_i(W_{p,n})}{\text{dist}(\lambda_i(W_{p,n}), I_k)} ) \alpha_{I_k} |I_k| + \frac{\delta^{k_0 }}{ k- k_0}
\end{split}
\end{equation*}
with probability at least $1-n^{-C_1-100}$.

Summing over all intervals for $k\ge 10$ (say), we have 
$$
\sum_{k=10}^{k_0} 100 \delta^{k-7}+ \sum _{k \ge k_0} \frac{\delta^{k_0}}{ k- k_0} \le \delta. 
$$

Using Riemann integration of the principal value integral, we obtain
\begin{equation}
\begin{split}
y \sum_{I_k} (1+\frac{\lambda_i(W_{p,n})}{\text{dist}(\lambda_i(W_{p,n}), I_k)} )\alpha_{I_k} |I_k| &= | p.v.\int_{a}^b y \frac{x \rho_{MP,y}(x)}{x-\lambda_i(W_{p,n})} \ dx | +o(1)
\end{split}
\end{equation}
where (see \cite{KWcov} for details)
\begin{equation}\label{pv}
p.v.\int_{a}^b y \frac{x \rho_{MP,y}(x)}{x-\lambda_i(W_{p,n})} \ dx =\begin{cases}
    \sqrt{y}+o(1), & \text{if } |\lambda_i(W_{p,n}) -a | =o(1); \\
   -\sqrt{y}+o(1), & \text{if } |\lambda_i(W_{p,n}) -b | =o(1).
  \end{cases}
\end{equation}
follows from the explicit formula for the Stieltjes transform and from residue calculus. 

Now for the rest of eigenvalues such that $|\lambda_i(W_{p,n})-\lambda_j(W_{p,n-1})|\le |I_0|+|I_1|+\ldots+|I_{10}|\le 4\eta/\delta^{80}$. By Theorem \ref{thm:LMPL} and Cauchy interlacing law, the number of eigenvalues is at most ${T_+-T_-} \le 8n\eta/\delta^{80} =8C K^2 \log n/\delta^{88}$ with probability at least $1-n^{-C_1-100}$ for constant $C>0$ sufficiently large. Thus 
$$\frac{\sqrt{T_+ - T_-}}{C^{1.5} K \sqrt{\log n}} \le \frac{1}{ \delta^{44} C} \le \delta$$
again by choosing $C$ sufficiently large. From Lemma \ref{lem:singularco}, by comparing (\ref{total}), (\ref{partial}) and (\ref{pv}), one can conclude with probability at least $1-n^{-C_1-10}$,
$|x| \le \frac{C_2 K^2 \log n}{\sqrt{n}}.$ The conclusion of Theorem \ref{thm:singularvector} follows from symmetry and union bounds.


\end{document}